\documentclass[11pt,reqno]{amsart}

\usepackage[T1]{fontenc}
\usepackage[utf8]{inputenc}
\usepackage{amsmath,amssymb,amsthm,mathtools}
\usepackage{mathrsfs}
\usepackage{enumitem}
\usepackage{booktabs}
\usepackage[expansion=false]{microtype}

\usepackage[letterpaper,
            textwidth=0.85\paperwidth,
            hcentering,
            top=1in,
            bottom=1in,
            headsep=18pt,
            footskip=32pt]{geometry}

\usepackage[colorlinks=true,linkcolor=blue,citecolor=blue,urlcolor=blue]{hyperref}

\usepackage{parskip}
\usepackage{ragged2e}
\RaggedRight
\usepackage{titlesec}
\titleformat{\section}
    {\normalfont\bfseries\Large}{\thesection}{0.75em}{}
\titleformat{\subsection}
    {\normalfont\bfseries\large}{\thesubsection}{0.75em}{}
\titleformat{\subsubsection}
    {\normalfont\bfseries}{\thesubsubsection}{0.75em}{}
\titlespacing*{\section}{0pt}{1.5\baselineskip}{0.75\baselineskip}
\titlespacing*{\subsection}{0pt}{1.25\baselineskip}{0.5\baselineskip}

\makeatletter
\def\l@section{\@tocline{1}{0pt}{1pc}{2.5pc}{}}
\makeatother

\usepackage{fancyhdr}
\newcommand{\papershorttitle}{Divisor Sums with Spacing Constraints}
\fancypagestyle{paperstyle}{%
    \fancyhf{}
    \fancyhead[C]{\normalfont\scshape\small\papershorttitle}
    \fancyfoot[C]{\normalfont\small\thepage}

}
\fancypagestyle{plain}{%
    \fancyhf{}
    \fancyhead[C]{\normalfont\scshape\small\papershorttitle}
    \fancyfoot[C]{\normalfont\small\thepage}

}
\fancypagestyle{firstpage}{%
    \fancyhf{}
    \fancyhead[C]{\normalfont\scshape\small\papershorttitle}
    \fancyfoot[C]{\normalfont\small\thepage}

}

\setlist[itemize]{leftmargin=2em}

\newcommand{\N}{\mathbb N}

\newcommand{\Pplus}{P^{+}}

\newcommand{\sqf}{\mathrm{sqf}}

\newcommand{\G}{\mathcal G}

\newcommand{\M}{\mathcal M}

\newcommand{\Nset}{\mathcal N}

\theoremstyle{plain}
\newtheorem{theorem}{Theorem}[section]
\newtheorem{proposition}[theorem]{Proposition}
\newtheorem{lemma}[theorem]{Lemma}

\theoremstyle{definition}

\numberwithin{equation}{section}

\title[Divisor sums with spacing constraints]{Divisor sums with spacing constraints}
\author{}
\date{22 October 2026}

\begin{document}

\begin{abstract}
For $x,y>0$ let $S(x,y)$ denote the largest value of $\sum_{n\in A}\tau(n)$
over subsets $A$ of $[1,x]\cap\N$ whose elements are pairwise separated by
more than $y$ (here $\tau$ is the divisor function). We determine the order
of magnitude, and in a subrange the precise asymptotic, of $S(x,y)$ when the
spacing parameter $y$ ranges over a polylogarithmic scale
$y=(\log x)^{1+\lambda(\log\lambda-1)}$, $\lambda>1$. We show that for every fixed $\lambda>2$ the spacing constraint costs a computable power of
$\log x$ in the exponent, while for every fixed $1<\lambda<2$ it costs
nothing at all: $S(x,y)\sim x\log x$, the same asymptotic as the
unconstrained sum $\sum_{n\le x}\tau(n)$. Near the threshold $\lambda=2$ the
two regimes meet, and we show that for any fixed $\kappa>0,$ $S(x,y)\asymp_{\kappa} x\log x$ uniformly as long
as $\lambda\le 2+2\kappa/\sqrt{\log\log x}$. The
proof combines the Sathe--Selberg distribution of the number of prime
factors, a smooth squarefree estimate obtained from a generating function
together with Cauchy's formula, a sieve estimate for correlations among the
squarefree parts of the integers that contribute most heavily to
$\sum_{n\le x}\tau(n)$, and two combinatorial devices (the Caro--Wei bound
and a concentration inequality on the symmetric group) used to extract a
well-spaced independent set from a sparse auxiliary graph.
\end{abstract}
\author{Aahan Chatterjee}
\address{D\'epartement de math\'ematiques et de statistique,
Universit\'e de Montr\'eal, Montr\'eal, Canada}
\maketitle
\tableofcontents
\section{Introduction}\label{sec:intro}

Dirichlet's classical divisor theorem states that
\begin{equation}\label{eq:dirichlet}
    \sum_{n\leq x}\tau(n)\sim x\log x.
\end{equation}
The question we are interested in is how much of this sum can remain if the integers that we select are required to
be some distance apart. More precisely, call a set $A\subseteq[1,x]\cap\mathbb N$
\emph{$y$-separated} if
\[
    |a-a'|>y
    \qquad\textit{for every }a,a'\in A,\ a\neq a'.
\]
Then we are interested in the size of the quantity:
\begin{equation}\label{eq:defS}
    S(x,y)
    :=
    \max_{\substack{A\subseteq[1,x]\cap\mathbb N\\
                     A\ \text{$y$-separated}}}
          \sum_{n\in A}\tau(n).
\end{equation} 
We study this question on a polylogarithmic scale.  Put
\[
    L:=\log\log x,
    \qquad
    I(t):=1+t(\log t-1),
\]
and, for $\lambda>1$, let
\begin{equation}\label{eq:scale}
    y_\lambda
    :=(\log x)^{I(\lambda)}
    =(\log x)^{1+\lambda(\log\lambda-1)}.
\end{equation}

The critical value $\lambda=2$ can be guessed from a simple
Sathe--Selberg heuristic.  From the anatomy of integers having large divisor contributions, we know that the divisor
sum is carried by integers having
\(\sim2\log\log x\) prime factors.  Sathe--Selberg and Stirling's formula suggest that
there are roughly
\[
    \approx\frac{x}{(\log x)^{2\log 2-1}}
\]
such integers (here and below, $\approx$ suppresses factors of size
$(\log x)^{o(1)}$, which do not affect the location of the threshold).  Thus these integers have average spacing
\[
    \approx(\log x)^{2\log 2-1}.
\]
If these integers were uniformly distributed in $[1,x]$, we would expect to be able to impose spacing below this scale without losing more
than a constant proportion of the divisor mass.  Above this scale, the
spacing condition should begin to force a genuine loss.  This predicts the
critical spacing
\[
    y=(\log x)^{2\log 2-1}=(\log x)^{I(2)}.
\]
Since $I(\lambda)$ is increasing for $\lambda>1$, the cases $\lambda<2$
and $\lambda>2$ correspond respectively to spacings below and above this
critical scale.  
This also explains why the transition occupies a window of length $1/\sqrt{L}$.  Repeating the computation with $2\log\log x$ replaced by
$(2+s)\log\log x$ gives $\approx x(\log x)^{-I(2+s)}$ integers, each carrying $\tau(n)\approx 2^{(2+s)L}$, and hence a total divisor weight of 
\[
    \approx x(\log x)^{(2+s)\log 2-I(2+s)}.
\]
At $s=0$ the exponent equals $2\log 2-I(2)=1$, in agreement with
\eqref{eq:dirichlet}, and by the discussion above $s=0$ is exactly where this
exponent is largest.  Taylor expanding $(2+s)\log 2-I(2+s)$ around $s=0$, we observe that the exponent is $1-cs^{2}+O(s^{3})$ for some $c>0$, so the divisor
weight becomes
\[
    \approx x\log x\cdot\exp\bigl(-cLs^{2}+O(Ls^{3})\bigr),
\]
The exponential factor changes by a constant amount when
$s$ is of order $L^{-1/2}$. Writing $\lambda=2+s$, this
suggests the transition window above $\lambda=2$ covered
by part~(ii) below.
Our main result establishes these heuristics rigorously.

\begin{theorem}\label{thm:main}
Let $L:=\log\log x$ and
\[
    y=y_\lambda:=(\log x)^{I(\lambda)},
    \qquad \lambda>1.
\]
\begin{enumerate}
    \item[\textup{(i)}]
    For every fixed $\lambda>2$,
    \[
        S(x,y)\asymp_\lambda
        \frac{x(\log x)^{\lambda\log 2}}
             {y(\log\log x)^{\log 2/(2\log\lambda)}}.
    \]

    \item[\textup{(ii)}]
    For every fixed $\kappa>0$, uniformly for
    \[
        1<\lambda\leq2+\frac{2\kappa}{\sqrt L},
    \]
    we have
    \[
        S(x,y)\asymp_\kappa x\log x.
    \]

    \item[\textup{(iii)}]
    For every fixed $1<\lambda<2$,
    \[
        S(x,y)\sim x\log x.
    \]
\end{enumerate}
\end{theorem}

Part~\textup{(ii)} includes the endpoint $\lambda=2$ and, in fact, an
$O(L^{-1/2})$ window above it.  There is no conflict between this and
part~\textup{(i)}, which is a fixed-$\lambda$ result and is not uniform as
$\lambda\downarrow2$.  Formally setting $\lambda=2$ in the expression in
part~\textup{(i)} gives $x\log x/\sqrt L$.  At the transition, however,
there are order $\sqrt L$ adjacent layers of this size, and their combined
contribution restores the factor $\sqrt L$.  The present paper does not give
a uniform result in the intermediate range
\[
    \lambda=2+\eta(x),
    \qquad
    \eta(x)\longrightarrow0,
    \qquad
    \eta(x)\sqrt L\longrightarrow\infty;
\]
for example, $\lambda=2+L^{-1/4}$ is not covered.  We also restrict to
$\lambda>1$, where $I$ is increasing.  At $\lambda=1$ the scale becomes
$y=1$, and a different formulation is more natural.

We are not aware of a previous treatment of the extremal quantity
$S(x,y)$.  There are, nevertheless, two nearby strands of work.  The
additive divisor problem studies shifted correlations such as
$\sum_{n\leq x}\tau(n)\tau(n+h)$, beginning with Ingham~\cite{Ingham} and
Estermann~\cite{Estermann}; see also~\cite{DI}.  Our problem asks for an
extremal packing rather than an average correlation, although the close-pair
estimate used below is related in spirit.  The second connection is with the
distribution of integers having divisors in prescribed intervals, especially
Ford's work~\cite{FordDiv}.  Our useful-divisor decomposition is motivated by
this point of view.

We briefly describe the proof.  The upper bounds in parts~\textup{(ii)} and
\textup{(iii)} follow immediately from \eqref{eq:dirichlet}.  For the upper
bound in part~\textup{(i)}, we split the admissible set at a suitably chosen
divisor threshold.  The contribution from integers below this threshold is
controlled by the cardinality bound $x/y+1$, while the contribution from
integers above it is estimated using the distribution of integers with many
prime factors.  More precisely, we write $n=m(n)r(n)$, where $m(n)$ is
powerful and $r(n)$ is squarefree.  The powerful part is handled by Rankin's
trick, while the squarefree upper tail is estimated using the uniform
Hardy--Ramanujan bound of Lemma~\ref{lem:hardy-ramanujan}.

For the lower bounds, we construct an auxiliary graph by placing an edge between two relevant integers whenever
their distance is at most $y$.  A set with the required separation is then an
independent set in this graph.  Proposition~\ref{prop:good-pair-count} shows
that the graph has few edges.  Its proof uses the structural result of
Proposition~\ref{prop:typical-useful-divisor}: for a positive proportion of
the squarefree integers under consideration---and, in the subcritical
argument, for all but a negligible proportion---the divisor $d_\ell(n)$ formed
from the prime factors at most $x^{1/\ell}$ satisfies
\[
    d_\ell(n)<x^{1/2-\delta}.
\]
After the two small divisors $d_1,d_2$ of a close pair are fixed, the
remaining $x^{1/\ell}$-rough cofactors lie in a one-parameter family.  The
bound $d_1d_2<x^{1-2\delta}$ makes the parameter interval sufficiently long
for an upper-bound sieve to be applied.  The exceptional set in the
useful-divisor result is controlled by the smooth squarefree estimate of
Lemma~\ref{lem:friable-ss}.

Once the close pairs have been counted, the Caro--Wei bound (see Theorem~\ref{cw}) gives the
independent sets needed in parts~\textup{(i)} and~\textup{(iii)}.  In the
transition range, it is not enough to find a large independent set in each
layer separately: the same independent set must have large components in every relevant
layer. To do this, we start with the probabilistic setup from the proof of
the Caro--Wei theorem: choose a uniformly random ordering $\sigma$ of the
vertices and let $A_v$ denote the event that $v$
precedes every other vertex of its closed neighbourhood
$N_G[v]$. Then
\[
    \mathcal I(\sigma):=\{v:A_v\text{ occurs}\}
\]
is an independent set. We apply a variant of the Azuma-Hoeffding inequality
on the symmetric group, stated in Theorem~\ref{ah} to each relevant
layer, followed by a union bound over these layers.
This will show that all the relevant components are simultaneously
large with positive probability.

Section~\ref{sec:standard} records the standard analytic estimates used in
the argument.  Sections~\ref{sec:friable} and~\ref{sec:useful-divisor} prove
the smooth estimate and the useful-divisor result, respectively.
Section~\ref{sec:correlation} establishes the close-pair estimate,
Section~\ref{sec:prelim} records the two combinatorial tools, and
Section~\ref{sec:main-proof} proves Theorem~\ref{thm:main}.

\section{Notation and Conventions}\label{sec:notation}

Throughout, \(x\ge 2\) denotes a real parameter. All asymptotic
statements are understood as \(x\to\infty\), unless an explicit range
of uniformity is stated. Whenever \(\sqrt L\) or \(\log L\) occurs,
\(x\) is understood to be sufficiently large. We write
\[
    L:=\log\log x,
    \qquad
    I(t):=1+t(\log t-1),
    \qquad
    y_t:=(\log x)^{I(t)}
\]
and usually take \(t>1\). Unless otherwise stated, a subscript on \(\ll\), \(\gg\), or \(\asymp\) records the permitted dependence of the implied constants. We write
\(\omega(n)\) for the number of distinct prime factors of \(n\) and
\(\tau(n)\) for the number of divisors of \(n\); an integer \(n\) is
\emph{squarefree} if \(\mu^2(n)=1\), and \emph{powerful} if \(p\mid n\)
implies \(p^2\mid n\).
We write \(\Pplus(n)\) for the largest prime factor of \(n\), with the convention \(\Pplus(1)=1\).
For \(X\ge1\) and \(k\ge0\), we write
\[
    N_k(X):=\{n\le X:\mu^2(n)=1,\ \omega(n)=k\},
\]
and use \(|N_k(X)|\) for its cardinality.

All graphs are finite, simple, and undirected. For a graph \(G\), we
write \(V(G)\) and \(E(G)\) for its vertex and edge sets. For
\(v\in V(G)\), let
\[
    N_G(v):=\{w\in V(G):\{v,w\}\in E(G)\},
    \qquad
    N_G[v]:=N_G(v)\cup\{v\},
\]
denote its open and closed neighbourhoods, respectively, and put
\[
    d_G(v):=|N_G(v)|.
\]
For \(W\subseteq V(G)\), \(G[W]\) denotes the subgraph induced by
\(W\). A set \(\mathcal I\subseteq V(G)\) is independent if no two
of its vertices are adjacent. We write \(\mathbb S_n\) for the
symmetric group on \(n\) elements; when considering a random
ordering of a graph's vertices, every ordering is equally likely.

For quick reference, Table~\ref{tab:notation} collects the recurring pieces
of notation used in the rest of the paper; each is (re)defined at its first
occurrence in the text.

\begin{table}[ht]
\centering
\renewcommand{\arraystretch}{1.35}
\begin{tabular}{@{}p{0.27\linewidth} p{0.69\linewidth}@{}}
\hline
Symbol & Meaning \\
\hline
\(L\) & \(\log\log x\) \\
\(\omega(n)\), \(\tau(n)\), \(\Pplus(n)\) & number of distinct prime factors, number of divisors, largest prime factor of \(n\) \\
\(\mu^2(n)=1\) & \(n\) is squarefree \\
\(m(n),\,r(n)\) & the powerful part and squarefree part of \(n=m(n)r(n)\), \eqref{eq:powerful-squarefree-decomp} \\
\(\Nset_{\kappa,M}(x)\) & \(n\le x\) with \(|\omega(r(n))-2L|\le\kappa\sqrt L\) and \(m(n)\le M\), \eqref{eq:def-N-kappa-M} \\
\(\Nset^{\sqf}_{\kappa}(x)\) & squarefree \(n\le x\) with \(|\omega(n)-2L|\le\kappa\sqrt L\), \eqref{eq:def-Nsf-kappa} \\
\(N_k(X)\) & squarefree \(n\le X\) with \(\omega(n)=k\) \\
\(F_k(x,y)\) & number of squarefree \(n\le x\) with \(\Pplus(n)\le y\) and \(\omega(n)=k\), \eqref{eq:def-Fkxy} \\
\(d_\ell(n)\) & the part of \(n\) built from its prime factors \(\le x^{1/\ell}\), \eqref{eq:def-dell} \\
\(\M_{k,\delta,\ell}(x)\) & \(n\in N_k(x)\) with \(d_\ell(n)<x^{1/2-\delta}\), Section~\ref{sec:useful-divisor} \\
\(\G_{\kappa,\delta,\ell}(x)\) & \(\bigcup_{|k-2L|\le\kappa\sqrt L}\M_{k,\delta,\ell}(x)\), Subsection~\ref{subsec:subcritical} \\
\(u\) & \(\log x/\log y\), for a smoothness parameter \(y\) \\
\(I(t)\) & \(1+t(\log t-1)\), the logarithmic scale function \\
\(y_t\) & \((\log x)^{I(t)}\) \\
\(\lambda\) & parameter governing the spacing scale \(y=y_\lambda\) \\
\(V(G)\), \(E(G)\) & vertex set and edge set of the graph \(G\) \\
\(d_G(v)\) & degree of \(v\) in \(G\) \\
\(N_G(v)\) & open neighbourhood of \(v\) in \(G\) \\
\(N_G[v]\) & closed neighbourhood \(N_G(v)\cup\{v\}\) \\
\(G[W]\) & subgraph of \(G\) induced by \(W\subseteq V(G)\) \\
\(\mathbb S_n\) & symmetric group on \(n\) elements \\
\(\mathcal I(\sigma)\) & vertices that precede every other vertex of their closed neighbourhood in the ordering \(\sigma\) \\
\hline
\end{tabular}
\caption{Recurring notation.}
\label{tab:notation}
\end{table}

Every integer \(n\) has a unique factorization
\begin{equation}\label{eq:powerful-squarefree-decomp}
    n=m(n)r(n),
\end{equation}
where \(m(n)\) is powerful, \(r(n)\) is squarefree, and \((m(n),r(n))=1\).  Equivalently,
\begin{equation*}
    m(n):=\prod_{\substack{p^\nu\Vert n\\ \nu\ge 2}}p^\nu,
    \qquad
    r(n):=\prod_{p\Vert n}p .
\end{equation*}
For fixed constants \(\kappa,M>0\), define
\begin{equation}\label{eq:def-N-kappa-M}
    \Nset_{\kappa,M}(x)
    :=
    \left\{n\le x:
    |\omega(r(n))-2L|\le \kappa\sqrt L,
    \ m(n)\le M
    \right\}.
\end{equation}
We also write
\begin{equation}\label{eq:def-Nsf-kappa}
    \Nset^{\sqf}_{\kappa}(x)
    :=
    \left\{n\le x:
    n\text{ squarefree},\ |\omega(n)-2L|\le \kappa\sqrt L
    \right\}.
\end{equation}

\section{Some standard inputs}\label{sec:standard}

This section collects some standard analytic estimates used later in the paper.
We first establish the explicit Mertens-product bounds needed for
the sieving arguments in Section~\ref{sec:correlation}. We then record upper bound and asymptotic estimates for squarefree integers with a prescribed number of prime factors. Finally, we prove the weighted concentration estimate used in the proof of part~\textup{(iii)} of Theorem~\ref{thm:main}.

\subsection{Explicit Mertens-product estimates}

The following interval estimate will be used in the verification of
Axiom~\textup{2} to enable the use of the fundamental lemma of sieve theory later in Section~\ref{sec:correlation}.

\begin{lemma}[Explicit Mertens-product bounds]
\label{lem:mertens-product-interval}
Let \(285<y_1\leq y_2\). Then
\begin{equation}\label{eq:mertens-product-interval}
    \prod_{y_1<p\leq y_2}
    \left(1-\frac1p\right)^{-1}
    \leq
    \left(1+\frac{2}{(\log y_1)^2}\right)
    \frac{\log y_2}{\log y_1}.
\end{equation}
Moreover, there is an absolute constant \(D_0>0\) such that
\begin{equation}\label{eq:global-mertens-upper}
    \prod_{p\leq y}
    \left(1-\frac1p\right)^{-1}
    \leq
    D_0\log y
\end{equation}
for every \(y\geq3/2\).
\end{lemma}

\begin{proof}
Let
\[
    M(y):=\prod_{p\leq y}\left(1-\frac1p\right),
\]
and let \(\gamma\) denote the Euler--Mascheroni constant. The explicit
estimates of Rosser and Schoenfeld
\cite[Theorem~7]{RosserSchoenfeld1962} give
\begin{equation}\label{eq:rosser-schoenfeld-mertens}
    \frac{e^{-\gamma}}{\log y}
    \left(1-\frac{1}{2(\log y)^2}\right)
    <
    M(y)
    <
    \frac{e^{-\gamma}}{\log y}
    \left(1+\frac{1}{2(\log y)^2}\right),
\end{equation}
where the upper bound is valid for \(y>1\), while the lower bound is
valid for \(y>285\).

For \(285<y_1\leq y_2\), we have
\[
    \prod_{y_1<p\leq y_2}
    \left(1-\frac1p\right)^{-1}
    =
    \frac{M(y_1)}{M(y_2)}.
\]
Applying the upper bound in
\eqref{eq:rosser-schoenfeld-mertens} to \(M(y_1)\) and the lower bound
to \(M(y_2)\), we obtain
\[
    \prod_{y_1<p\leq y_2}
    \left(1-\frac1p\right)^{-1}
    \leq
    \frac{\log y_2}{\log y_1}
    \left(\frac{
        1+\dfrac{1}{2(\log y_1)^2}
    }{
        1-\dfrac{1}{2(\log y_2)^2}
    }\right).
\]
Set
\[
    a:=\frac{1}{2(\log y_1)^2},
    \qquad
    b:=\frac{1}{2(\log y_2)^2}.
\]
Then
\[
    0\leq b\leq a<\frac12.
\]
Consequently,
\[
\begin{aligned}
    \frac{1+a}{1-b}
    &=
    1+\frac{a+b}{1-b}\\
    &\leq
    1+\frac{2a}{1-a}\\
    &\leq
    1+4a\\
    &=
    1+\frac{2}{(\log y_1)^2}.
\end{aligned}
\]
This proves \eqref{eq:mertens-product-interval}.

It remains to prove \eqref{eq:global-mertens-upper}. If \(y>285\),
the lower bound in \eqref{eq:rosser-schoenfeld-mertens} gives
\[
\begin{aligned}
    \prod_{p\leq y}
    \left(1-\frac1p\right)^{-1}
    &<
    e^\gamma\log y
    \left(1-\frac{1}{2(\log y)^2}\right)^{-1}\\
    &\leq
    e^\gamma
    \left(1-\frac{1}{2(\log 285)^2}\right)^{-1}
    \log y.
\end{aligned}
\]
For \(3/2\leq y\leq285\), we have
\[
    \prod_{p\leq y}
    \left(1-\frac1p\right)^{-1}
    \leq
    \prod_{p\leq285}
    \left(1-\frac1p\right)^{-1},
\]
while
\[
    \log y\geq\log(3/2)>0.
\]
Thus \eqref{eq:global-mertens-upper} follows, for example, upon taking
\[
    D_0
    :=
    \max\left\{
        e^\gamma
        \left(1-\frac{1}{2(\log 285)^2}\right)^{-1},
        \frac{1}{\log(3/2)}
        \prod_{p\leq285}
        \left(1-\frac1p\right)^{-1}
    \right\}.
\]
\end{proof}

\subsection{Squarefree integers with a prescribed number of prime factors}

Recall that
\[
    N_k(X)
    :=
    \{n\leq X:\mu^2(n)=1,\ \omega(n)=k\}.
\]
We shall use two estimates for \(\lvert N_k(X)\rvert\). The first is
a uniform upper bound valid for every \(k\). The second is an
asymptotic estimate in the Sathe--Selberg range
\(k\ll\log\log X\).

\begin{lemma}[Hardy--Ramanujan upper bound]
\label{lem:hardy-ramanujan}
There exist absolute constants \(c_1,c_2>0\) such that
\begin{equation}\label{eq:HR-uniform}
    \lvert N_k(X)\rvert
    \leq
    c_1\frac{X}{\log X}
    \frac{(\log\log X+c_2)^{k-1}}{(k-1)!}
\end{equation}
for all \(X\geq2\) and all integers \(k\geq1\).
\end{lemma}

\begin{proof}
Estimate \cite[(11.2)]{GranvilleANT} proves
\eqref{eq:HR-uniform} for the larger set
\[
    \{n\leq X:\omega(n)=k\}.
\]
Since \(N_k(X)\) is the squarefree subset of this set, the same upper
bound holds for \(N_k(X)\).
\end{proof}

The Hardy--Ramanujan estimate will be used when uniformity over all
\(k\) is needed, particularly when controlling the tails of the
distribution. For \(k\) in a bounded multiple of \(\log\log X\), we
use the sharper Sathe--Selberg asymptotic.

For real \(t\geq0\), define
\begin{equation}\label{eq:def-sathe-G}
    G(t)
    :=
    \frac{1}{\Gamma(t+1)}
    \prod_p
    \left(1+\frac{t}{p}\right)
    \left(1-\frac1p\right)^t.
\end{equation}
For every fixed \(R>0\),
\[
    \log\left(
        \left(1+\frac{t}{p}\right)
        \left(1-\frac1p\right)^t
    \right)
    =
    O_R(p^{-2})
\]
uniformly for \(0\leq t\leq R\). Since \(\sum_p p^{-2}<\infty\),
the series of logarithms converges absolutely and uniformly on
\([0,R]\). Each Euler factor and \(\Gamma(t+1)\) is positive there,
so \(G\) is continuous and strictly positive. In particular, there
are constants \(c_R,C_R>0\) such that
\[
    c_R\leq G(t)\leq C_R
    \qquad(0\leq t\leq R).
\]

\begin{lemma}[Sathe--Selberg]\label{sathe_uniform}
For every fixed \(R>0\), there exists \(X_0=X_0(R)\geq2\) such that,
uniformly for \(X\geq X_0\) and integers \(k\) satisfying
\[
    1\leq k\leq R\log\log X,
\]
we have
\begin{equation}\label{eq:sathe_uniform}
    \lvert N_k(X)\rvert
    =
    \frac{X}{\log X}
    \frac{(\log\log X)^{k-1}}{(k-1)!}
    \left\{
        G\left(\frac{k-1}{\log\log X}\right)
        +
        O_R\left(
            \frac{k}{(\log\log X)^2}
        \right)
    \right\}.
\end{equation}
\end{lemma}

\begin{proof}
This is the first specialization listed after
\cite[Theorem~4, pp.~86--87]{SelbergSathe}. The corresponding Euler
product is
\[
    \prod_p(1+zp^{-s})
    =
    \sum_{n\geq1}\frac{\mu^2(n)z^{\omega(n)}}{n^s}
    \qquad(\Re s>1).
\]
Thus its coefficient of \(z^k\), summed over \(n\leq X\), counts
the squarefree set \(N_k(X)\). Selberg allows an arbitrary fixed
range constant for this specialization.
\end{proof}

Since \(G\) is bounded away from zero on \([0,R]\), the additive
error in \eqref{eq:sathe_uniform} can equivalently be written as a
relative error \(1+O_R(k/(\log\log X)^2)\). We record separately the
two-sided consequence used later.

\begin{lemma}
\label{lem:sathe-standard}
For every fixed \(B>0\), there exists \(X_0=X_0(B)\geq2\) such that
\begin{equation}\label{eq:sathe-standard}
    \lvert N_k(X)\rvert
    \asymp_B
    \frac{X}{\log X}
    \frac{(\log\log X)^{k-1}}{(k-1)!}
\end{equation}
uniformly for \(X\geq X_0\) and integers \(k\) satisfying
\[
    1\leq k\leq B\log\log X.
\]
\end{lemma}

\begin{proof}
Put \(L_X:=\log\log X\) and apply Lemma~\ref{sathe_uniform} with
\(R=B\). For \(1\leq k\leq BL_X\), we have
\[
    0\leq\frac{k-1}{L_X}\leq B,
    \qquad
    \frac{k}{L_X^2}\leq\frac{B}{L_X}.
\]
The bounds on \(G\) established above therefore apply. Thus the
expression in braces in \eqref{eq:sathe_uniform} is bounded above
and below by positive constants depending only on \(B\), provided
\(X_0(B)\) is sufficiently large. This proves
\eqref{eq:sathe-standard}.
\end{proof}

\subsection{Explicit Stirling Inequality}
Now we also record an explicit form of Stirling's inequalities (see \cite{Robbins}):
\begin{lemma}\label{lem:robbins}
    For every integer $n\geq1$,
\begin{equation}\label{eq:robbins-stirling}
    \sqrt{2\pi}\,n^{n+1/2}e^{-n}e^{1/(12n+1)}
    <
    n!
    <
    \sqrt{2\pi}\,n^{n+1/2}e^{-n}e^{1/(12n)}.
\end{equation}
\end{lemma}
\subsection{Weighted concentration}

The final input controls the total divisor weight of integers lying
outside the central range used in the construction of the good set.
We derive the required estimate from the classical variance estimate for $\omega(n)$ and the convergence of the divisor-weighted sum over powerful integers.

\begin{lemma}[Weighted concentration]
\label{lem:weighted-central}
For every \(\varepsilon>0\), there exist constants
\[
    \kappa=\kappa(\varepsilon)\geq1
    \qquad\text{and}\qquad
    M=M(\varepsilon)\geq1
\]
such that, for all sufficiently large \(x\),
\begin{equation}\label{eq:weight_central}
    \sum_{\substack{n\leq x\\
                    n\notin\Nset_{\kappa,M}(x)}}
    \tau(n)
    \leq
    \varepsilon x\log x.
\end{equation}
\end{lemma}

\begin{proof}
The sum
\[
    C_{\mathrm{pow}}
    :=
    \sum_{\substack{m\geq1\\m\ \mathrm{powerful}}}
       \frac{\tau(m)}m
    =
    \prod_p\left(
       1+\sum_{\nu\geq2}\frac{\nu+1}{p^\nu}
    \right)
    <\infty
\]
converges, since each local factor is \(1+O(p^{-2})\).
Using the canonical factorization \(n=m(n)r(n)\), and then
discarding the coprimality and squarefreeness restrictions for an
upper bound, we obtain
\[
\begin{aligned}
    \sum_{\substack{n\leq x\\m(n)>M}}\tau(n)
    &\leq
    \sum_{\substack{M<m\leq x\\m\ \mathrm{powerful}}}
       \tau(m)\sum_{r\leq x/m}\tau(r)\\
    &\ll
    x\log x
    \sum_{\substack{m>M\\m\ \mathrm{powerful}}}
       \frac{\tau(m)}m.
\end{aligned}
\]
Here we used \(\sum_{r\leq T}\tau(r)\ll T\log(2T)\) for \(T\geq1\).
Choose \(M=M(\varepsilon)\geq1\) so that this contribution is at
most \(\varepsilon x\log x/2\).

We next establish the second-moment estimate
\begin{equation}\label{eq:weighted-squarefree-second-moment}
    \sum_{r\leq T}
       \mu^2(r)2^{\omega(r)}(\omega(r)-2L)^2
    \ll T\log x\,L
\end{equation}
uniformly for \(x^{1/2}\leq T\leq x\), with an absolute implied
constant, once \(x\) is sufficiently large.
Put \(H:=\log\log T\) and
\[
    A_T(u):=\sum_{a\leq u}(\omega(a)-H)^2.
\]
The mean of \(\omega(a)\) for \(a\leq u\) is
\(\log\log u+O(1)\). From
\cite[(9.2)]{GranvilleANT}, we obtain
\[
    \sum_{a\leq u}(\omega(a)-\log\log u)^2
    \ll u\log\log(3u).
\]
This holds for \(u\geq3\), after enlarging the absolute implied
constant if necessary.
Using \((s+t)^2\leq2s^2+2t^2\), we can change the center in the second moment bound as follows:
\[
    A_T(u)
    \ll u\left\{
       \log\log(3u)
       +\left(\log\frac{\log T}{\log u}\right)^2
    \right\}
    \qquad(3\leq u\leq T).
\]
Partial summation therefore gives
\[
\begin{aligned}
    \sum_{a\leq T}\frac{(\omega(a)-H)^2}{a}
    &=\frac{A_T(T)}{T}
       +\int_3^T\frac{A_T(u)}{u^2}\,du
       +O(H^2)\\
    &\ll \log T\,H
       +\int_3^T
          \left(\log\frac{\log T}{\log u}\right)^2
          \frac{du}{u}\\
    &\ll \log T\,H.
\end{aligned}
\]
Here the lower-endpoint terms are \(O(H^2)\), which is
\(O(H\log T)\) for large \(T\), and, on setting
\(v:=\log u/\log T\),
\[
    \int_3^T
       \left(\log\frac{\log T}{\log u}\right)^2
       \frac{du}{u}
    \leq \log T\int_0^1\left(\log\frac1v\right)^2\,dv
    =O(\log T).
\]
Since \(L-\log2\leq H\leq L\) in the stated range, again using $(s+t)^2\le 2s^2+2t^2,$ we have
\[
\begin{aligned}
    \sum_{a\leq T}\frac{(\omega(a)-L)^2}{a}
    &\leq
       2\sum_{a\leq T}\frac{(\omega(a)-H)^2}{a}
       +2(H-L)^2\sum_{a\leq T}\frac1a\\
    &\ll \log T\,(H+1)\ll \log x\,L.
\end{aligned}
\]

For squarefree \(r\), each ordered factorization \(r=ab\)
satisfies \(\omega(r)=\omega(a)+\omega(b)\), and there are
\(2^{\omega(r)}\) such factorizations. Hence, by symmetry and
\((u+v)^2\leq2u^2+2v^2\),
\[
\begin{aligned}
    &\sum_{r\leq T}
       \mu^2(r)2^{\omega(r)}(\omega(r)-2L)^2\\
    &\quad=
       \sum_{ab\leq T}\mu^2(ab)
          (\omega(a)+\omega(b)-2L)^2\\
    &\quad\leq
       4\sum_{ab\leq T}(\omega(a)-L)^2\\
    &\quad\leq
       4T\sum_{a\leq T}\frac{(\omega(a)-L)^2}{a}
       \ll T\log x\,L.
\end{aligned}
\]
This proves \eqref{eq:weighted-squarefree-second-moment}.

For \(x\geq M^2\), every \(m\leq M\) satisfies
\(x/m\geq x^{1/2}\). Since \(\tau(r)=2^{\omega(r)}\) for squarefree
\(r\), applying \eqref{eq:weighted-squarefree-second-moment} gives
\[
\begin{aligned}
    &\sum_{\substack{n\leq x\\m(n)\leq M\\
          |\omega(r(n))-2L|>\kappa\sqrt L}}\tau(n)\\
    &\qquad\leq
    \frac1{\kappa^2L}
    \sum_{\substack{m\leq M\\m\ \mathrm{powerful}}}
       \tau(m)
       \sum_{r\leq x/m}
          \mu^2(r)2^{\omega(r)}(\omega(r)-2L)^2\\
    &\qquad\ll
    \frac{x\log x}{\kappa^2}
    \sum_{\substack{m\leq M\\m\ \mathrm{powerful}}}
       \frac{\tau(m)}m
    \ll
    \frac{C_{\mathrm{pow}}}{\kappa^2}\,x\log x.
\end{aligned}
\]
Choose \(\kappa=\kappa(\varepsilon)\geq1\) sufficiently large that
this contribution is at most \(\varepsilon x\log x/2\).
Combining the two bounds proves \eqref{eq:weight_central}.
\end{proof}
\section{A Smooth squarefree estimate}\label{sec:friable}

This section proves a uniform upper bound for the number of \(y\)-smooth
squarefree integers up to \(x\) with exactly \(k\) prime factors, in the
range \(k=O(\log\log y)\). Recall that the unrestricted counting function
of \(y\)-smooth integers is governed, in the classical range, by
\[
    \Psi(x,y)\sim x\rho(u),
    \qquad
    u:=\frac{\log x}{\log y},
\]
where \(\rho\) denotes the Dickman--de Bruijn function (see \cite[Chapter III.5]{Ten}). This suggests heuristically that imposing the condition \(P^{+}(n)\le y\) in addition to the Sathe--Selberg condition $\omega(n)=k$ would introduce a Dickman-type saving and give an upper bound of the shape:
\[
    F_k(x,y)
    \stackrel{?}{\ll}
    \rho(u)\,
    \frac{x}{\log x}
    \frac{(\log\log y)^k}{k!}.
\]
However, for our purposes, it is enough to obtain the weaker exponential saving
\(e^{-Au}\), for arbitrary fixed \(A>0\). We prove such a bound directly
from a generating-function estimate using a Rankin-type argument (see \cite[Chapter~III.5]{Ten}) and Cauchy's formula.

For \(x\ge y\ge 2\) and \(k\ge 0\), set
\begin{equation}\label{eq:def-Fkxy}
    F_k(x,y)
    :=
    \#\{n\le x:\mu^2(n)=1,\ \Pplus(n)\le y,\ \omega(n)=k\}.
\end{equation}
The estimate below is the technical input used to show that most relevant squarefree integers possess a useful small divisor.

\begin{lemma}\label{lem:friable-ss}
Let \(B>0\) and \(A>0\) be fixed. There is a constant
\(C=C(A,B)>0\) such that, whenever \(x>y\ge C\), and
\[
    u:=\frac{\log x}{\log y},
\]
one has, uniformly for positive integers
\[
    1\le k\le B\log\log y,
\]
\begin{equation}\label{eq:friable-ss-bound}
    F_k(x,y)
    \ll_{A,B}
    \frac{\sqrt k\,x e^{-Au}}{\log x}
    \frac{(\log\log y)^k}{k!}.
\end{equation}
\end{lemma}

We first record the generating-function estimate from which Lemma~\ref{lem:friable-ss} will follow.

\begin{lemma}\label{lem:Sz-bound}
Let \(B>0\) and \(A>0\) be fixed and denote
\[
    S_z(x,y):=
    \sum_{\substack{n\le x\\ \Pplus(n)\le y}}z^{\omega(n)}\mu^2(n).
\]
There exists a constant \(C=C(A,B)>0\) such that for \(x>y\ge C\), one has, uniformly for \(|z|\le B\),
\begin{equation}\label{eq:Sz-bound}
    |S_z(x,y)|
    \ll_{A,B}
    \frac{x}{u}e^{-Au}(\log y)^{|z|-1},
    \qquad
    u=\frac{\log x}{\log y}.
\end{equation}
\end{lemma}

To prove Lemma~\ref{lem:Sz-bound}, we use the following standard
mean-value estimate for non-negative multiplicative functions (see \cite[Corollary~III.3.6]{Ten}).

\begin{proposition}\label{prop:tenenbaum-mean}
Let \(f\) be a non-negative multiplicative function satisfying
\[
    0\le f(p^\nu)\le \beta_1\beta_2^{\nu-1},
    \qquad \nu\ge1,
\]
for fixed constants \(\beta_1>0\) and \(0\le\beta_2<2\).
Then, uniformly for \(X\ge2\),
\begin{equation}\label{eq:tenenbaum-mean}
    \sum_{n\le X}f(n)
    \ll_{\beta_1,\beta_2}
    X
    \prod_{p\le X}
    \left(1-\frac1p\right)
    \left(
        \sum_{\nu\ge0}\frac{f(p^\nu)}{p^\nu}
    \right).
\end{equation}
\end{proposition}
We are now ready to prove Lemma~\ref{lem:Sz-bound}.
\begin{proof}[Proof of Lemma~\ref{lem:Sz-bound}]
Choose \(C=C(A,B)\) sufficiently large that
\begin{equation}\label{eq:y-lower-A}
    \log y\ge 4A
\end{equation}
whenever \(y\ge C\).

If \(z=0\), only the term \(n=1\) contributes, so \(S_0(x,y)=1\). On the other hand,
\[
    \frac{x}{u}e^{-Au}(\log y)^{-1}
    =
    \frac{x}{\log x}e^{-Au}.
\]
Since \(\log y\ge 4A\), we have
\[
    Au\le \frac14\log x,
\]
and hence
\[
    \frac{x}{\log x}e^{-Au}
    \ge \frac{x^{3/4}}{\log x}>1
\]
after enlarging \(C(A,B)\) if necessary. Thus the result is immediate
when \(z=0\), and we may henceforth assume \(z\neq0\).

We first dispose of the range \(y>x^{3/4}\). In this case
\[
    1<u=\frac{\log x}{\log y}<\frac43,
    \qquad
    \log y\asymp\log x.
\]
Apply Proposition~\ref{prop:tenenbaum-mean} with \(X=x\) to the
non-negative multiplicative function
\[
    f(n)
    :=
    |z|^{\omega(n)}\mu^2(n)
    \mathbf 1_{\{\Pplus(n)\le y\}}.
\]
Since \(f(p)\le B\) and \(f(p^\nu)=0\) for \(\nu\ge2\), the hypotheses
of Proposition~\ref{prop:tenenbaum-mean} hold, for instance with
\(\beta_1=B\) and \(\beta_2=1\). Hence
\begin{align*}
    |S_z(x,y)|
    &\le
    \sum_{\substack{n\le x\\ \Pplus(n)\le y}}
    |z|^{\omega(n)}\mu^2(n) \\
    &\ll_B
    x
    \prod_{y<p\le x}\left(1-\frac1p\right)
    \prod_{p\le y}
    \left(1-\frac1p\right)
    \left(1+\frac{|z|}{p}\right).
\end{align*}
For \(p\le y\),
\[
    \left(1-\frac1p\right)
    \left(1+\frac{|z|}{p}\right)
    =
    1+\frac{|z|-1}{p}
    -\frac{|z|}{p^2}\le 1+\frac{|z|-1}{p}
\]
Therefore, using \(1+t\le e^t\) for \(t\in \mathbb{R}\),
\[
\begin{aligned}
    \prod_{p\le y}
    \left(1-\frac1p\right)
    \left(1+\frac{|z|}{p}\right)
    &\le
    \exp\left(
        (|z|-1)\sum_{p\le y}\frac1p
    \right) \\
    &\ll_B
    (\log y)^{|z|-1},
\end{aligned}
\]
where in the last step we used Mertens' estimate
\[
    \sum_{p\le y}\frac1p
    =
    \log\log y+O(1).
\]
Also, by Mertens' product theorem,
\[
    \prod_{y<p\le x}\left(1-\frac1p\right)
    \asymp
    \frac{\log y}{\log x}
    =
    \frac1u.
\]
Consequently,
\[
    |S_z(x,y)|
    \ll_B
    \frac{x}{u}(\log y)^{|z|-1}.
\]
Since \(1<u<4/3\),
\[
    e^{-Au}\ge e^{-4A/3},
\]
and hence
\[
    1\le e^{4A/3}e^{-Au}.
\]
Thus, after enlarging the implied constant by a factor depending only on
\(A\),
\[
    |S_z(x,y)|
    \ll_{A,B}
    \frac{x}{u}e^{-Au}(\log y)^{|z|-1}.
\]
This proves the desired estimate when \(y>x^{3/4}\).

We may therefore assume for the remainder of the proof that
\[
    y\le x^{3/4}.
\]
Set
\[
    \alpha:=\frac{4A}{3\log y}.
\]
By Rankin's trick,
\[
    \mathbf 1_{\{n>x^{3/4}\}}
    \le
    \left(\frac{n}{x^{3/4}}\right)^\alpha
    =
    x^{-3\alpha/4}n^\alpha.
\]
Hence
\begin{equation}\label{eq:Sz-decomp}
\begin{aligned}
    |S_z(x,y)|
    &\le
    \sum_{\substack{n\le x^{3/4}\\ \Pplus(n)\le y}}
    |z|^{\omega(n)}\mu^2(n) \\
    &\qquad+
    x^{-3\alpha/4}
    \sum_{\substack{n\le x\\ \Pplus(n)\le y}}
    n^\alpha |z|^{\omega(n)}\mu^2(n).
\end{aligned}
\end{equation}

We first estimate the second term in \eqref{eq:Sz-decomp}. Apply
Proposition~\ref{prop:tenenbaum-mean}, with \(X=x\), to
\[
    f(n)
    :=
    n^\alpha |z|^{\omega(n)}\mu^2(n)
    \mathbf 1_{\{\Pplus(n)\le y\}}.
\]
For \(p\le y\), we have
\[
    p^\alpha
    \le
    y^\alpha
    =
    e^{4A/3}.
\]
Thus
\[
    f(p)\le Be^{4A/3},
    \qquad
    f(p^\nu)=0\quad(\nu\ge2),
\]
so the hypotheses of Proposition~\ref{prop:tenenbaum-mean} hold with $\beta_1=Be^{4A/3},\beta_2=1.$
We obtain
\begin{align}
    &\sum_{\substack{n\le x\\ \Pplus(n)\le y}}
    n^\alpha |z|^{\omega(n)}\mu^2(n)
    \notag\\
    &\qquad\ll_{A,B}
    x
    \prod_{y<p\le x}\left(1-\frac1p\right)
    \prod_{p\le y}
    \left(1-\frac1p\right)
    \left(1+|z|p^{\alpha-1}\right).
    \label{eq:weighted-tenenbaum}
\end{align}

We now estimate the Euler product over \(p\le y\). For such \(p\),
\begin{align*}
    \left(1-\frac1p\right)
    \left(1+|z|p^{\alpha-1}\right)
    &=
    1-\frac1p
    +|z|p^{\alpha-1}
    -|z|p^{\alpha-2} \\
    &=
    1+\frac{|z|-1}{p}
    +|z|\frac{p^\alpha-1}{p}
    -|z|p^{\alpha-2}.
\end{align*}
Since
\[
    0\le\alpha\log p\le\frac{4A}{3},
\]
the mean-value theorem applied to the exponential function gives
\[
    p^\alpha-1
    =
    e^{\alpha\log p}-1
    \ll_A
    \alpha\log p
    \ll_A
    \frac{\log p}{\log y}.
\]
Moreover \(p^\alpha\ll_A1\), so
\[
    p^{\alpha-2}\ll_A\frac1{p^2}.
\]
It follows that
\[
    \left(1-\frac1p\right)
    \left(1+|z|p^{\alpha-1}\right)
    =
    1+\frac{|z|-1}{p}+E_p,
\]
where
\[
    |E_p|
    \ll_{A,B}
    \frac{\log p}{p\log y}+\frac1{p^2}.
\]
Since \(1+t\le e^t\) for every \(t\in\mathbb R\), we obtain
\begin{align*}
    &\prod_{p\le y}
    \left(1-\frac1p\right)
    \left(1+|z|p^{\alpha-1}\right)\\
    &\qquad\le
    \exp\left(
        \sum_{p\le y}
        \left(
            \frac{|z|-1}{p}+E_p
        \right)
    \right)\\
    &\qquad\le
    \exp\left(
        (|z|-1)\sum_{p\le y}\frac1p
        +
        O_{A,B}\left(
            \frac1{\log y}
            \sum_{p\le y}\frac{\log p}{p}
            +
            \sum_{p\le y}\frac1{p^2}
        \right)
    \right).
\end{align*}
By Mertens' estimates,
\[
    \sum_{p\le y}\frac1p
    =
    \log\log y+O(1),
    \qquad
    \sum_{p\le y}\frac{\log p}{p}
    =
    \log y+O(1),
\]
while
\[
    \sum_p\frac1{p^2}<\infty.
\]
Consequently,
\begin{equation}\label{eq:small-prime-product}
    \prod_{p\le y}
    \left(1-\frac1p\right)
    \left(1+|z|p^{\alpha-1}\right)
    \ll_{A,B}
    (\log y)^{|z|-1}.
\end{equation}

Again by Mertens' product theorem,
\begin{equation}\label{eq:large-prime-product}
    \prod_{y<p\le x}\left(1-\frac1p\right)
    \asymp
    \frac{\log y}{\log x}
    =
    \frac1u.
\end{equation}
Combining
\eqref{eq:weighted-tenenbaum},
\eqref{eq:small-prime-product}, and
\eqref{eq:large-prime-product}, we obtain
\[
    \sum_{\substack{n\le x\\ \Pplus(n)\le y}}
    n^\alpha |z|^{\omega(n)}\mu^2(n)
    \ll_{A,B}
    \frac{x}{u}(\log y)^{|z|-1}.
\]
Since
\[
    x^{-3\alpha/4}
    =
    \exp\left(
        -\frac{3\alpha}{4}\log x
    \right)
    =
    \exp\left(
        -A\frac{\log x}{\log y}
    \right)
    =
    e^{-Au},
\]
the second term in \eqref{eq:Sz-decomp} is
\begin{equation}\label{eq:second-Sz-term}
    \ll_{A,B}
    \frac{x}{u}e^{-Au}(\log y)^{|z|-1}.
\end{equation}

It remains to estimate the first term in \eqref{eq:Sz-decomp}. Apply
Proposition~\ref{prop:tenenbaum-mean}, now with \(X=x^{3/4}\), to
\[
    f(n)
    :=
    |z|^{\omega(n)}\mu^2(n)
    \mathbf 1_{\{\Pplus(n)\le y\}}.
\]
Since \(y\le x^{3/4}\), Mertens' product theorem gives
\[
    \prod_{y<p\le x^{3/4}}
    \left(1-\frac1p\right)
    \asymp
    \frac{\log y}{\log(x^{3/4})}
    \asymp
    \frac1u.
\]
The same calculation as in \eqref{eq:small-prime-product}, now with \(\alpha=0\), gives
\[
    \prod_{p\le y}
    \left(1-\frac1p\right)
    \left(1+\frac{|z|}{p}\right)
    \ll_B
    (\log y)^{|z|-1}.
\]
Thus
\begin{equation}\label{eq:first-Sz-term-pre}
    \sum_{\substack{n\le x^{3/4}\\ \Pplus(n)\le y}}
    |z|^{\omega(n)}\mu^2(n)
    \ll_B
    \frac{x^{3/4}}{u}(\log y)^{|z|-1}.
\end{equation}
Finally, by \eqref{eq:y-lower-A},
\[
    Au
    =
    A\frac{\log x}{\log y}
    \le
    \frac14\log x,
\]
and hence
\[
    x^{3/4}\le xe^{-Au}.
\]
It follows from \eqref{eq:first-Sz-term-pre} that the first term in
\eqref{eq:Sz-decomp} also satisfies
\[
    \ll_{A,B}
    \frac{x}{u}e^{-Au}(\log y)^{|z|-1}.
\]
Combining this with \eqref{eq:second-Sz-term} proves
\eqref{eq:Sz-bound}.
\end{proof}

\begin{proof}[Proof of Lemma~\ref{lem:friable-ss}]
Let \(C=C(A,B)\) be as in Lemma~\ref{lem:Sz-bound}, and suppose that
\(x>y\ge C\).
Let us write
\[
    L_y:=\log\log y,
    \qquad
    R:=\frac{k}{L_y}.
\]
Since \(k\ge1\), we have \(R>0\), while the hypothesis
\(k\le B L_y\) gives
\[
    0<R\le B.
\]
In particular, Lemma~\ref{lem:Sz-bound} applies on the circle
\(|z|=R\). 

Since
\[
    S_z(x,y)=\sum_{j\ge0}F_j(x,y)z^j,
\]
Cauchy's formula gives
\[
    F_k(x,y)
    =
    \frac{1}{2\pi i}
    \int_{|z|=R}
    \frac{S_z(x,y)}{z^{k+1}}\,dz.
\]
Hence, by Lemma~\ref{lem:Sz-bound},
\begin{equation}\label{eq:F*k}
    \begin{split}
    F_k(x,y)
    &\ll_{A,B}
    \frac{x}{u}e^{-Au}
    (\log y)^{R-1}R^{-k}\\
    &=
    \frac{x}{\log x}e^{-Au}
    (\log y)^R R^{-k},
    \end{split}
\end{equation}
since \(u\log y=\log x\).
Now
\[
    (\log y)^R
    =
    e^{R\log\log y}
    =
    e^k,
\]
and therefore
\[
    (\log y)^R R^{-k}
    =
    e^k\left(\frac{L_y}{k}\right)^k.
\]
Now using \eqref{eq:robbins-stirling}:
\[
    (\log y)^R R^{-k}=e^k\left(\frac{L_y}{k}\right)^k
    \ll
    \sqrt{k}\frac{L_y^k}{k!}.
\]
So substituting this into \eqref{eq:F*k} proves \eqref{eq:friable-ss-bound}.
\end{proof}
\section{Typical squarefree integers have a useful divisor}\label{sec:useful-divisor}

The correlation estimates of Section~\ref{sec:correlation} will be proved by
writing each integer under consideration as a product of a \emph{small} divisor \(d_{\ell}(n)\), built from its prime factors below \(x^{1/\ell}\), and a cofactor free of prime factors below \(x^{1/\ell}\). The condition
\[
    d_\ell(n)<x^{1/2-\delta}
\]
is what makes the correlation estimates of the next section possible:
after fixing the small divisors associated with a close pair, the rough cofactors lie in a one-parameter family whose parameter ranges over an
interval sufficiently long compared with the sifting level
\(x^{1/\ell}\).  This makes it possible to apply an upper-bound sieve to the
cofactors.

The purpose of this section is therefore to show that the useful-divisor
condition holds for a sufficiently large proportion of the squarefree
integers under consideration.  Part~\textup{(a)} shows that
\(\ell\) need not grow with \(x\): once \(\ell\) is chosen to be a
sufficiently large fixed constant (depending on the relevant parameters),
a positive proportion of \(N_k(x)\) is already good.  This fixed-\(\ell\)
statement will be particularly useful in the graph-theoretic arguments in
parts~\textup{(i)} and~\textup{(ii)} of the main theorem.  Part~\textup{(b)}
gives a stronger quantitative refinement when \(\ell\) is at least a
sufficiently large multiple of \(\log k\). In this case, the exceptional proportion is exponentially small in \(\ell\). This decay will be useful to prove ~\textup{(iii)} of the main theorem.
We refer to the integers satisfying
\(d_\ell(n)<x^{1/2-\delta}\) as \emph{good}.

Now recall that
\[
    N_k(x)=\{n\le x:\mu^2(n)=1,\ \omega(n)=k\}.
\]
For every \(n\in N_k(x)\), define
\begin{equation}\label{eq:def-dell}
    d_{\ell}(n)
    :=
    \prod_{\substack{p\mid n\\ p\le x^{1/\ell}}}p.
\end{equation}
Fix \(\delta\in(0,1/2)\). For integers \(1\le \ell<k\), let
\(\M_{k,\delta,\ell}(x)\) be the set of \(n\in N_k(x)\) for which
\[
    d_{\ell}(n)<x^{1/2-\delta}.
\]
By definition, every prime divisor of \(d_{\ell}(n)\) is at most
\(x^{1/\ell}\), while every prime divisor of \(n/d_{\ell}(n)\) is greater than
\(x^{1/\ell}\).

\begin{proposition}\label{prop:typical-useful-divisor}
Let \(\gamma>1\) and \(\delta\in(0,1/2)\) be fixed, let
\(L:=\log\log x\), and suppose that
\[
    L<k<\gamma L.
\]
\begin{enumerate}[label=\textup{(\alph*)}]
    \item Uniformly for \(1\le \ell<k\),
          \begin{equation}\label{eq:pos-prop-divisor}
    |N_k(x)\setminus\M_{k,\delta,\ell}(x)|
    \ll_{\gamma,\delta}
    \frac{|N_k(x)|}{\ell}.
\end{equation}
In particular, if \(\ell\) is chosen sufficiently large in terms of
\(\gamma\) and \(\delta\), then
\[
    |\M_{k,\delta,\ell}(x)|\gg_{\gamma,\delta}|N_k(x)|.
\]
 
    \item For every \(A>0,\) uniformly for
    \[
        \log k\le \ell<k,
    \]
    we have that
    \begin{equation}\label{eq:useful-divisor-quant}
        |N_k(x)\setminus \M_{k,\delta,\ell}(x)|
        \ll_{A,\delta,\gamma}
        \sqrt{k}\,e^{-A\ell}|N_k(x)|.
    \end{equation}
    In particular, if \(A>1/2\), for such \(\ell\),
    \[
        |\M_{k,\delta,\ell}(x)|=(1-o(1))|N_k(x)|,
    \]
    where the \(o(1)\) is as \(k\to\infty\).
\end{enumerate}
\end{proposition}
 
\begin{proof}
We prove the two assertions separately.
 
\begin{enumerate}[label=\textup{(\alph*)}]
    \item We set \(y:=x^{1/\ell}\). For \(\ell=1\),
\[
    |N_k(x)\setminus\M_{k,\delta,1}(x)|
    \le |N_k(x)|,
\]
so the desired estimate is immediate. We may therefore assume
\(\ell\ge2\). Since $n\in N_k(x)$ is square-free, we have the inequality
    \[
        \#\{n\in N_k(x):p\mid n\}
        \le |N_{k-1}(x/p)|.
    \]
    Since \(p\le y\le x^{1/2}\), we have \(x/p\ge x^{1/2}\), and so,
uniformly for such \(p\),
\[
    \log\log(x/p)=L+O(1).
\]
Thus using this in conjunction with the inequality
\[1+t\le e^t\] for every $t\in \mathbb{R},$ we obtain:
\begin{equation}\label{eq:alizwell}
    \begin{split}
        (\log\log(x/p)+O(1))^{k-2}&=(L+O(1))^{k-2}\\
        =L^{k-2}(1+O(1/L))^{k-2}&\ll_{\gamma} L^{k-2}
    \end{split}
\end{equation} where for the last inequality we have used the inequality
\begin{equation*}
    1+t \le e^t
\end{equation*} for every $t\in \mathbb{R}$ and the fact that $k-2<k\le \gamma L.$
Now recalling that $x/p\ge x^{1/2}\ge 2,$ for $x\ge 4$, we can use \eqref{eq:HR-uniform} for $X=x/p.$ This gives:
    \begin{align*}
        |N_{k-1}(x/p)|
        &\ll
        \frac{x/p}{\log(x/p)}
        \frac{(\log\log(x/p)+O(1))^{k-2}}{(k-2)!} \\
        &\ll_{\gamma}
        \frac{x/p}{\log(x/p)}
        \frac{L^{k-2}}{(k-2)!} \\
        &\le
        \frac1p\frac{2x}{\log x}
        \frac{L^{k-1}}{(k-1)!}\frac{k-1}{L}.
    \end{align*}
    For all sufficiently large \(x\) in terms of \(\gamma\), and
\(1\le k<\gamma L\),
\eqref{eq:sathe-standard} applies with \(k\), \(X=x\) and
\(B=\gamma\), yielding:
\begin{equation*}
    |N_k(x)|\asymp_{\gamma}\frac{x}{\log{x}}\frac{(\log\log{x})^{k-1}}{(k-1)!}
\end{equation*} and hence
\begin{equation*}
    \frac1p\frac{2x}{\log x}
        \frac{L^{k-1}}{(k-1)!}\frac{k-1}{L}\ll_{\gamma} \frac1p |N_k(x)|
\end{equation*}
where the last step again uses $(k-1)/L<\gamma.$
    Hence
    \begin{align*}
        \sum_{n\in N_k(x)}\log d_{\ell}(n)
        &=
        \sum_{p\le y}
        \#\{n\in N_k(x):p\mid n\}\log p \\
        &\ll_{\gamma}
        |N_k(x)|\sum_{p\le y}\frac{\log p}{p}.
    \end{align*}
    By Mertens' estimate,
\[
    \sum_{p\le y}\frac{\log p}{p}
    =\log y+O(1),
\]
and therefore
\[
    \sum_{n\in N_k(x)}\log d_\ell(n)
    \ll_\gamma
    |N_k(x)|\log y
    =
    |N_k(x)|\frac{\log x}{\ell}.
\]
Now, if
\(n\in N_k(x)\setminus\M_{k,\delta,\ell}(x)\), then
\[
    \log d_\ell(n)
    \ge
    \left(\frac12-\delta\right)\log x.
\]
Hence Markov's inequality gives
\[
\begin{aligned}
    \left(\frac12-\delta\right)\log x\,
    |N_k(x)\setminus\M_{k,\delta,\ell}(x)|
    &\le
    \sum_{n\in N_k(x)}\log d_\ell(n)\\
    &\ll_\gamma
    |N_k(x)|\frac{\log x}{\ell}.
\end{aligned}
\]
Consequently,
\[
    |N_k(x)\setminus\M_{k,\delta,\ell}(x)|
    \ll_{\gamma,\delta}
    \frac{|N_k(x)|}{\ell},
\]
which proves \eqref{eq:pos-prop-divisor}.
 
    \item Let us again set \(y:=x^{1/\ell}\) and suppose $A>0.$ Now write
    \[
        n=d_{\ell}(n)r_n.
    \]
    Then every prime factor of \(r_n\) is greater than \(y\). If
    \(j:=\omega(r_n)\ge1\), then
    \[
        r_n>y^j,
        \qquad
        d_{\ell}(n)=\frac{n}{r_n}<\frac{x}{y^j}=x^{1-j/\ell}.
    \]
    Moreover \(j<\ell\), since otherwise \(r_n>y^\ell=x\). If \(j=0\), then
    \(r_n=1\) and \(d_{\ell}(n)=n\). We know that $\omega(d_{\ell}(n))=k-\omega(r_n).$ For a fixed \(1\le j\le \ell-1\), the number of exceptional \(n\) with
\(\omega(r_n)=j\) is at most the number of pairs \((d,r)\) satisfying
\[
    x^{1/2-\delta}\le d<x^{1-j/\ell},
    \qquad
    \mu^2(d)=1,
    \qquad
    \omega(d)=k-j,
    \qquad
    \Pplus(d)\le y,
\]
together with
\[
    r\le \frac{x}{d},
    \qquad
    \left(r,\prod_{p\le y}p\right)=1.
\]
Here we have dropped the additional condition
\[
    \omega(r)=j,
\]
which can only enlarge the count. If $j=0$, the condition becomes
$x^{1/2-\delta}\leq n\leq x$. Thus the total number of exceptional
integers is at most the sum of the following two quantities:
    \begin{equation}\label{eq:exceptional-first-term}
        \sum_{j=1}^{\ell-1}
        \sum_{\substack{
            x^{1-j/\ell}>d\ge x^{1/2-\delta}\\
            \omega(d)=k-j\\
            \Pplus(d)\le y
        }}
        \mu^2(d)
        \sum_{\substack{
            r\le x/d\\
            (r,\prod_{p\le y}p)=1
        }}1,
    \end{equation}
    and the second is
    \begin{equation}\label{eq:exceptional-second-term}
        \sum_{\substack{
            x\ge d\ge x^{1/2-\delta}\\
            \omega(d)=k\\
            \Pplus(d)\le y
        }}
        \mu^2(d).
    \end{equation}
 
    In \eqref{eq:exceptional-first-term}, we have \(d<x^{1-j/\ell}\le
x^{1-1/\ell}=x/y\), and hence \(x/d>y=x^{1/\ell}>1\) for $x>1.$ So we can use Fan's inequality for \(y\)-rough integers (see
\cite{FanRough}):
\[
    \#\left\{r\le X:
        \left(r,\prod_{p\le y}p\right)=1
      \right\}
    \le \frac{X}{\log y}
\] uniformly for $X\ge y>1$
to obtain
\[
        \sum_{\substack{
            r\le x/d\\
            (r,\prod_{p\le y}p)=1
        }}1
        \le
        \frac{x}{d\log y}
        =
        \frac{x\ell}{d\log x}
\] uniformly for all $d<x^{1-1/\ell},\quad x>1.$
Let us denote
    \[
        a:=x^{1/2-\delta},
        \qquad
        b_j:=x^{1-j/\ell},
    \]
    and let \(N\) be the largest integer \(1\le N\le\ell-1\) for which
    \(b_N>a\). For \(j>N\) the interval \(a\le d<b_j\) is empty, so only
    \(1\le j\le N\) contributes.
    Thus \eqref{eq:exceptional-first-term} is
    \begin{equation}\label{eq:first-after-sieve}
        \ll\sum_{j=1}^{N}\frac{x\ell}{\log x}\sum_{\substack{b_j>d\ge a\\\omega(d)=k-j\\\Pplus(d)\le y}}\frac{\mu^2(d)}{d}.
    \end{equation}

    For fixed \(j\), partial summation gives the upper bound
    \begin{align}\label{eq:partial-summation-Sj}
        \sum_{\substack{
            a\le d<b_j\\
            \omega(d)=k-j\\
            \Pplus(d)\le y
        }}
        \frac{\mu^2(d)}{d}
        &\le
        \frac{F_{k-j}(b_j,y)}{b_j}
        +
        \int_a^{b_j}\frac{F_{k-j}(t,y)}{t^2}\,dt.
    \end{align}
 
    Let
\[
    L_y:=\log\log y=L-\log\ell.
\]
Since \(\ell<k<\gamma L\), we have
\[
    L_y=L+O_{\gamma}(\log L),
\]
and hence, for \(x\) sufficiently large in terms of $\gamma$,
\[
    1\le k-j\le k<\gamma L\le2\gamma L_y.
\]
Moreover, the hypothesis \(\ell\ge \log k\), together with
\(k>L\), implies
\[
    y=x^{1/\ell}<x^{1/2-\delta}=a\le t
\]
for all \(x\) sufficiently large in terms of $\delta$, uniformly for \(t\in[a,b_j]\).
Finally let \(A_1>0\) be a parameter to be chosen later and let
\[
    C_{\mathrm{fr}}=C_{\mathrm{fr}}(A_1,2\gamma)
\]
be the constant appearing in Lemma~\ref{lem:friable-ss}.
Since
\[
    y=x^{1/\ell}
    \ge x^{1/k}
    \ge x^{1/(\gamma L)}
    \longrightarrow\infty,
\]
we have \(y\ge C_{\mathrm{fr}}\) for all \(x\) sufficiently large in terms of $A_1,\gamma$. Therefore for all $x$ sufficiently large in terms of $A_1,\gamma,\delta,$ Lemma~\ref{lem:friable-ss} applies uniformly to \(F_{k-j}(t,y)\), with \(A=A_1\) and \(B=2\gamma\), and yields
    \begin{align*}
        \frac{x\ell}{\log x}
        \int_a^{b_j}\frac{F_{k-j}(t,y)}{t^2}\,dt
        &\ll_{A_1,\gamma}
        \frac{x\ell\sqrt{k}}{\log x}
        \frac{(\log\log y)^{k-j}}{(k-j)!}
        \\
        &\qquad\times
        \int_a^{b_j}
        \frac{\exp(-A_1\log t/\log y)}{t\log t}\,dt.
    \end{align*}
    With the change of variables
    \[
        v:=\frac{\log t}{\log y},
    \]
    the last integral is bounded by
    \[
    \begin{aligned}
        \int_{(1/2-\delta)\ell}^{\log b_j/\log y}
        \frac{e^{-A_1v}}{v}\,dv
        &\le
        \frac{1}{(1/2-\delta)\ell}
        \int_{(1/2-\delta)\ell}^{\infty}e^{-A_1v}\,dv\\
        &\ll_{A_1,\delta}
        \frac{e^{-A_1(1/2-\delta)\ell}}{\ell}.
    \end{aligned}
    \]
    Hence
    \begin{equation}\label{eq:integral-Sj-bound}
        \frac{x\ell}{\log x}
        \int_a^{b_j}\frac{F_{k-j}(t,y)}{t^2}\,dt
        \ll_{A_1,\gamma,\delta}
        \frac{x\sqrt{k}}{\log x}
        e^{-A_1(1/2-\delta)\ell}
        \frac{(\log\log y)^{k-j}}{(k-j)!}.
    \end{equation}
     We also observe that Lemma~\ref{lem:friable-ss} at \(t=b_j\) gives
    \begin{align*}
        \frac{x\ell}{\log x}\frac{F_{k-j}(b_j,y)}{b_j}
        &\ll_{A_1,\gamma}
        \frac{x\ell\sqrt{k}}{\log x\,\log b_j}
        e^{-A_1\log b_j/\log y}
        \frac{(\log\log y)^{k-j}}{(k-j)!} \\
        &\ll_{A_1,\gamma,\delta}
        \frac{x\sqrt{k}}{\log x}
        e^{-A_1(1/2-\delta)\ell}
        \frac{(\log\log y)^{k-j}}{(k-j)!},
    \end{align*}
    since \(\log b_j\ge(1/2-\delta)\log x\) and \(\ell\ll L=o(\log x)\). Hence the endpoint term satisfies the same bound.
 
    We now compare the remaining sum with \(|N_k(x)|\). Let us denote
    \(L_y:=\log\log y\). Since \(L_y\le L\) and $k<\gamma L,$ by
\eqref{eq:sathe-standard} applied with $B=\gamma$, we have that:
    \begin{equation*}
        |N_k(x)|\asymp_{\gamma}\frac{x}{\log{x}}\cdot\frac{L^{k-1}}{(k-1)!}
    \end{equation*}
    and hence uniformly for \(1\le j\le N\):
    \begin{align*}
        \frac{x}{\log x}\frac{L_y^{k-j}}{(k-j)!}
        &\ll_{\gamma}
        |N_k(x)|
        \frac{(k-1)!}{(k-j)!L^{j-1}}
        \ll_{\gamma} \gamma^{j-1}|N_k(x)|
    \end{align*} since
    \begin{equation*}
        \begin{split}
            \frac{(k-1)!}{(k-j)!L^{j-1}}&=\prod_{r=k-j+1}^{k-1}\left(\frac{r}{L}\right)\le \gamma^{j-1}
        \end{split}
    \end{equation*} since $r \le k-1<k\le \gamma L.$
    Consequently,
    \begin{equation}\label{eq:sum-over-j-ratio}
        \sum_{j=1}^{N}
        \frac{x}{\log x}
        \frac{(\log\log y)^{k-j}}{(k-j)!}
        \ll_\gamma
        |N_k(x)|\sum_{j=1}^{N}\gamma^{j-1}
        \ll_\gamma
        \gamma^{\ell}|N_k(x)|.
    \end{equation}
    Therefore \eqref{eq:first-after-sieve} is
    \begin{equation}\label{eq:first-exceptional-final}
        \ll_{A_1,\gamma,\delta}
        \sqrt{k}\,\gamma^{\ell}
        e^{-A_1(1/2-\delta)\ell}|N_k(x)|.
    \end{equation}
    Let us choose \(A_1>0\) sufficiently large so that
\[
    \left(\frac12-\delta\right)A_1>A+\log\gamma.
\]
Then
\[
    \gamma^\ell
    e^{-A_1(1/2-\delta)\ell}
    \le e^{-A\ell},
\]
so \eqref{eq:first-exceptional-final} is
\[
    \ll_{A,\delta,\gamma}
    \sqrt{k}\,e^{-A\ell}|N_k(x)|.
\]
 
    By definition, the second term \eqref{eq:exceptional-second-term} is at most \(F_k(x,y)\). Now again recall that $1\le k \le 2\gamma L_y$ and hence Lemma~\ref{lem:friable-ss} applies with $B=2\gamma.$ So applying for a parameter \(A_2>0\), gives
    \[
        F_k(x,y)
        \ll_{A_2,\gamma}
        \frac{\sqrt{k}x e^{-A_2\ell}}{\log x}
        \frac{(\log\log y)^k}{k!}.
    \]
    Since \(k<\gamma L,\) applying \eqref{eq:sathe-standard} with $B=\gamma$ and using $k>L$ gives
    \[
        \frac{x}{\log x}\frac{(\log\log y)^k}{k!}
        \ll_{\gamma}
        \frac{L}{k}|N_k(x)|
        <|N_k(x)|.
    \]
    Hence
    \[
        F_k(x,y)
        \ll_{A_2,\gamma}
        \sqrt{k}e^{-A_2\ell}|N_k(x)|.
    \]
    Taking \(A_2>A\) proves \eqref{eq:useful-divisor-quant}.
    Finally, if \(A>1/2\), then the hypothesis \(\ell\ge\log k\) gives
\[
    \sqrt{k}\,e^{-A\ell}
    \le k^{1/2-A}
    =o(1).
\]
Hence
\[
    |\M_{k,\delta,\ell}(x)|
    =(1-o(1))|N_k(x)|,
\]
which proves the final assertion.
\end{enumerate}
\end{proof}

\section{A correlation estimate for the good squarefree sets}\label{sec:correlation}

This section bounds the number of pairs
\[
    (m,n)\in
    \M_{i,\delta,\ell}(x)\times\M_{j,\delta,\ell}(x)
\]
with \(0<|am-bn|\leq H\), where \(a,b\) are bounded. These bounds
control the edges of the auxiliary graphs in
Section~\ref{sec:main-proof}.

The proof writes each good integer as its small divisor times an
\(x^{1/\ell}\)-rough cofactor. After fixing the two small divisors
\(d_1,d_2\), the linear relation between the cofactors gives a
one-parameter family. The condition
\[
    d_1d_2<x^{1-2\delta}
\]
makes the parameter interval sufficiently long to apply the
fundamental lemma of sieve theory.

\begin{proposition}\label{prop:good-pair-count}
Fix a real number \(K\geq1\), \(\delta\in(0,1/2)\), and
\(\gamma>1\). There exist an integer
\[
    \ell_0=\ell_0(\delta)\geq1
\]
and a constant \(C_\gamma>0\), depending only on \(\gamma\), such
that, for all sufficiently large \(x\) in terms of the fixed
parameters $K,\gamma,\delta$, uniformly for integers \(a,b,i,j,\ell\) and real \(H\)
satisfying
\[
    1\leq a,b\leq K,
    \qquad
    L<i,j<\gamma L,
    \qquad
    H\geq1,
    \qquad
    \ell_0\leq\ell\leq L,
\]
one has
\begin{equation}\label{eq:single-pair-count}
\begin{aligned}
    &\#\left\{
        (m,n)\in
        \M_{i,\delta,\ell}(x)\times\M_{j,\delta,\ell}(x):
        0<|am-bn|\leq H,\ bn\leq x
    \right\}\\
    &\qquad\ll_{K,\delta,\gamma}
    \frac{C_\gamma^\ell\ell^2H}{abx}
    |N_i(x)|\,|N_j(x)|.
\end{aligned}
\end{equation}
\end{proposition}
 
\subsection{Upper-bound sieve for two linear forms}
Before proceeding with the proof, we will prove a sieving bound for quadratic polynomials which can be written as a product of two linear forms. This will be a key input for the proof of Proposition~\ref{prop:good-pair-count}.

\begin{lemma}[Upper-bound sieve for two linear forms]
\label{lem:two-linear-form-sieve}
 Suppose \(I\subset\mathbb R\) is an interval containing
\(X:=\#(I\cap\mathbb Z)\) integers, and let
\[
    f(t)=(\alpha_1t+\beta_1)(\alpha_2t+\beta_2),
    \qquad
    \alpha_i,\beta_i\in\mathbb Z.
\]
For \(p\) prime, define
\[
    \nu_f(p)
    :=
    \#\{t\bmod p:f(t)\equiv0\pmod p\}.
\]
Then there exist absolute constants \(c_1,c_2>0\) such that, uniformly
for all \(X\geq c_1\) and \(2\leq Y\leq X\) satisfying
\[
    \log X\geq c_2\log Y,
\]
\begin{equation}\label{eq:two-linear-form-sieve}
    \#\left\{t\in I\cap\mathbb Z:\left(f(t),\prod_{p\le Y}p\right)=1\right\}
    \le
    5X\prod_{p\le Y}
    \left(1-\frac{\nu_f(p)}p\right).
\end{equation}
\end{lemma}

We prove Lemma~\ref{lem:two-linear-form-sieve} using the Fundamental
Lemma of Sieve Theory in the form given by
\cite[Theorem~18.11\textup{(b)}]{KoukoulopoulosPrimes}, specialized to
sieving by all primes \(p\leq Y\). We first recall the relevant setup
and hypotheses.

Let \(\mathcal A=(a_n)_{n\in\mathcal I}\) be a finite sequence of integers. For each positive integer \(d\), define
\[
    A_d:=\#\{n\in\mathcal I:d\mid a_n\}.
\]
The corresponding sifted sum is
\[
    S(\mathcal A,Y)
:=
\#\{n\in\mathcal I:
      p\nmid a_n\text{ for every prime }p\leq Y\}.
\]
Throughout, products indexed by \(p\) are taken over primes.

The Fundamental Lemma is stated in terms of the following three
axioms.

\medskip
\noindent
\textbf{Axiom 1.}
There exist a parameter \(X>1\), a nonnegative multiplicative function
\(\nu\), and real numbers \(r_d\) such that
\begin{equation}\label{eq:sieve-axiom-one}
    A_d=\frac{\nu(d)}{d}X+r_d
\end{equation}
for every squarefree integer \(d\) all of whose prime factors are at
most \(Y\). Moreover,
\[
    \nu(p)<p
    \qquad
    \text{for every prime }p\leq Y.
\]

\medskip
\noindent
\textbf{Axiom 2.}
There exist constants \(\kappa\geq0\) and \(C>0\) such that
\begin{equation}\label{eq:sieve-axiom-two}
    \prod_{y_1<p\leq y_2}
    \left(1-\frac{\nu(p)}{p}\right)^{-1}
    \leq
    \left(1+\frac{C}{\log y_1}\right)
    \left(\frac{\log y_2}{\log y_1}\right)^\kappa
\end{equation}
uniformly for
\[
    \frac32\leq y_1\leq y_2\leq Y.
\]

\medskip
\noindent
\textbf{Axiom 3.}
There exist constants \(A>0\) and \(m\in\mathbb N\), together with a
parameter \(D\geq1\), such that
\begin{equation}\label{eq:sieve-axiom-three}
    \sum_{\substack{d\leq D\\ d\ \mathrm{squarefree}\\
        p\mid d\,\Longrightarrow\,p\leq Y
        \text{ for every prime }p}}
    \tau_m(d)|r_d|
    \leq
    \frac{X}{(\log X)^A}.
\end{equation}
Here
\[
    \tau_m(n)
    :=
    \#\{(d_1,\ldots,d_m)\in\mathbb N^m:
          d_1\cdots d_m=n\}
\]
is the \(m\)-fold divisor function. In particular,
\(\tau_1(n)=1\).

We now state the part of the Fundamental Lemma that will be used.

\begin{theorem}[Fundamental Lemma of Sieve Theory]
\label{thm:fundamental-sieve}
Suppose that the sequence \(\mathcal A\) satisfies
Axioms~\textup{1} and~\textup{2} for some fixed
\(\kappa,C>0\). Define
\[
    u_\kappa
    :=
    1+\frac{2}{e^{0.53/\kappa}-1}.
\]
Suppose, in addition, that Axiom~\textup{3} holds with
\[
    m=1,
    \qquad
    A=\kappa+1,
    \qquad
    D\geq Y^{u_\kappa}.
\]
If
\[
    \log X\gg\log Y
\]
and \(X\) and \(D\) are sufficiently large in terms of \(\kappa\) and
\(C\), then
\begin{equation}\label{eq:fundamental-sieve-comparison}
    S(\mathcal A,Y)
    \leq
    5X
    \prod_{p\leq Y}
    \left(1-\frac{\nu(p)}{p}\right).
\end{equation}
\end{theorem}

\begin{proof}[Proof of Lemma~\ref{lem:two-linear-form-sieve}]
Set
\[
    \mathcal A:=\bigl(f(t)\bigr)_{t\in I\cap\mathbb Z}.
\]
Suppose first that
\[
    \nu_f(p)=p
\]
for some prime \(p\leq Y\). Then \(p\mid f(t)\) for every
\(t\in\mathbb Z\), so
\[
    S(\mathcal A,Y)=0.
\]
The corresponding factor in the Euler product is also zero. Therefore
\eqref{eq:two-linear-form-sieve} holds trivially.

We may henceforth assume that
\[
    \nu_f(p)<p
    \qquad
    \text{for every prime }p\leq Y.
\]
We verify the hypotheses of
Theorem~\ref{thm:fundamental-sieve} with sieve dimension
\(\kappa=2\).

We recall that
\[
    \nu_f(d)
    :=
    \#\{t\bmod d:f(t)\equiv0\pmod d\}
\]
and hence by the Chinese remainder theorem, \(\nu_f\) is multiplicative.
In particular, if \(d\) is squarefree, then
\[
    \nu_f(d)=\prod_{p\mid d}\nu_f(p).
\]

Let
\[
    \mathcal S_d
    :=
    \{r\bmod d:f(r)\equiv0\pmod d\}.
\]
Then
\[
    d\mid f(t)
    \quad\Longleftrightarrow\quad
    t\equiv r\pmod d
    \quad\text{for some }r\in\mathcal S_d.
\]
Since \(I\cap\mathbb Z\) consists of \(X\) consecutive integers, for
each \(r\in\mathcal S_d\) we have
\[
    \#\{t\in I\cap\mathbb Z:t\equiv r\pmod d\}
    =
    \frac{X}{d}+O(1).
\]
Consequently,
\[
\begin{aligned}
    A_d
    &:=
    \#\{t\in I\cap\mathbb Z:d\mid f(t)\} \\
    &=
    \sum_{r\in\mathcal S_d}
    \#\{t\in I\cap\mathbb Z:t\equiv r\pmod d\} \\
    &=
    \frac{\nu_f(d)}{d}X+r_d,
    \qquad
    |r_d|\ll\nu_f(d).
\end{aligned}
\]
This verifies Axiom~\textup{1}.

Under our assumption that \(\nu_f(p)<p\), neither linear factor of
\(f\) is identically zero modulo \(p\). Since a nonzero linear
polynomial has at most one root modulo \(p\), it follows that
\[
    0\leq\nu_f(p)\leq2
    \qquad
    \text{for every prime }p\leq Y.
\]
We verify Axiom~\textup{2} separately according as \(y_1>285\) or
\(y_1\leq285\).

Suppose first that
\[
    285<y_1\leq y_2\leq Y.
\]
For \(p>2\), we have
\[
    \left(1-\frac{\nu_f(p)}p\right)^{-1}
    \leq
    \left(1-\frac2p\right)^{-1},
\]
and
\[
    \left(1-\frac2p\right)^{-1}
    =
    \left(1-\frac1p\right)^{-2}
    \left(1+\frac{1}{p(p-2)}\right).
\]
Consequently,
\begin{align}\label{eq:axiom-two-large-y1}
    &\prod_{y_1<p\leq y_2}
    \left(1-\frac{\nu_f(p)}p\right)^{-1} \notag\\
    &\qquad\leq
    \left(
        \prod_{y_1<p\leq y_2}
        \left(1-\frac1p\right)^{-1}
    \right)^2
    \prod_{y_1<p\leq y_2}
    \left(1+\frac{1}{p(p-2)}\right).
\end{align}

We next estimate the second product. Since \(1+t\leq e^t\) for
\(t\geq0\), we have
\[
\begin{aligned}
    1
    &\leq
    \prod_{y_1<p\leq y_2}
    \left(1+\frac{1}{p(p-2)}\right)\\
    &\leq
    \exp\left(
        \sum_{p>y_1}\frac{1}{p(p-2)}
    \right).
\end{aligned}
\]
Since \(y_1>285\), we have \(p-2\geq p/2\) for every prime
\(p>y_1\), and therefore
\[
    \sum_{p>y_1}\frac{1}{p(p-2)}
    \leq
    2\sum_{n>y_1}\frac1{n^2}
    \ll
    \frac1{y_1}.
\]
Using the Taylor expansion of the exponential function at the origin,
we conclude that
\begin{equation}\label{eq:convergent-local-correction}
    \prod_{y_1<p\leq y_2}
    \left(1+\frac{1}{p(p-2)}\right)
    =
    1+O\left(\frac1{y_1}\right),
\end{equation}
uniformly in \(y_2\).

By Lemma~\ref{lem:mertens-product-interval},
\[
    \prod_{y_1<p\leq y_2}
    \left(1-\frac1p\right)^{-1}
    \leq
    \left(1+\frac{2}{(\log y_1)^2}\right)
    \frac{\log y_2}{\log y_1}.
\]
Substituting this and
\eqref{eq:convergent-local-correction} into
\eqref{eq:axiom-two-large-y1}, we obtain
\[
\begin{aligned}
    \prod_{y_1<p\leq y_2}
    \left(1-\frac{\nu_f(p)}p\right)^{-1}
    &\leq
    \left(1+\frac{2}{(\log y_1)^2}\right)^2
    \left(1+O\left(\frac1{y_1}\right)\right)
    \left(\frac{\log y_2}{\log y_1}\right)^2\\
    &\leq
    \left(1+\frac{C_0}{\log y_1}\right)
    \left(\frac{\log y_2}{\log y_1}\right)^2
\end{aligned}
\]
for some absolute constant \(C_0>0\).

It remains to consider
\[
    \frac32\leq y_1\leq285,
    \qquad
    y_1\leq y_2\leq Y.
\]
The possible factor corresponding to \(p=2\) is at most \(2\). For
\(p>2\), the preceding factorization gives
\[
\begin{aligned}
    \prod_{y_1<p\leq y_2}
    \left(1-\frac{\nu_f(p)}p\right)^{-1}
    &\leq
    2
    \prod_{\substack{p\leq y_2\\p>2}}
    \left(1-\frac2p\right)^{-1}\\
    &\leq
    2
    \left(
        \prod_{p\leq y_2}
        \left(1-\frac1p\right)^{-1}
    \right)^2
    \prod_{p>2}
    \left(1+\frac{1}{p(p-2)}\right).
\end{aligned}
\]
The final Euler product converges, and
Lemma~\ref{lem:mertens-product-interval} gives
\[
    \prod_{p\leq y_2}
    \left(1-\frac1p\right)^{-1}
    \ll\log y_2
\]
because \(y_2\geq y_1\geq3/2\). Hence there is an absolute constant \(D_1>0\) such that
\begin{equation}\label{eq:axiom-two-small-y1}
    \prod_{y_1<p\leq y_2}
    \left(1-\frac{\nu_f(p)}p\right)^{-1}
    \leq
    D_1(\log y_2)^2.
\end{equation}

Choose \(C\geq C_0\) sufficiently large that
\[
    \frac{1+C/\log 285}{(\log 285)^2}
    \geq D_1.
\]
The function
\[
    u\longmapsto\frac{1+C/u}{u^2}
\]
is decreasing for \(u>0\). Therefore, when
\(3/2\leq y_1\leq285\),
\[
    D_1
    \leq
    \frac{1+C/\log 285}{(\log 285)^2}
    \leq
    \frac{1+C/\log y_1}{(\log y_1)^2}.
\]
Combining this with \eqref{eq:axiom-two-small-y1}, we obtain
\[
    \prod_{y_1<p\leq y_2}
    \left(1-\frac{\nu_f(p)}p\right)^{-1}
    \leq
    \left(1+\frac{C}{\log y_1}\right)
    \left(\frac{\log y_2}{\log y_1}\right)^2.
\]
Thus Axiom~\textup{2} holds with sieve dimension \(\kappa=2\).
Finally, choose
\[
    D:=\frac{X}{(\log X)^5}.
\]
Since 
\[
    \nu_f(d)
    =
    \prod_{p\mid d}\nu_f(p)
    \leq
    2^{\omega(d)}=\tau(d)
\]
for squarefree \(d\), we obtain
\[
\begin{aligned}
    \sum_{\substack{d\leq D\\ d\ \mathrm{squarefree}\\
        p\mid d\,\Longrightarrow\,p\leq Y
        \text{ for every prime }p}}
    |r_d|
    &\ll
    \sum_{d\leq D}\mu^2(d)\tau(d) \\
    &\ll
    D\log{D} \\
    &\ll
    \frac{X}{(\log X)^4}
    \leq
    \frac{X}{(\log X)^3}
\end{aligned}
\]
for all sufficiently large \(X\). Hence Axiom~\textup{3} holds with
\[
    m=1,
    \qquad
    A=3=\kappa+1.
\]

Choose the absolute constant \(c_2>u_2\) sufficiently large to meet the
condition \(\log X\gg\log Y\) in
Theorem~\ref{thm:fundamental-sieve}. Since
\[
    \log X\geq c_2\log Y,
\]
we have \(Y\leq X^{1/c_2}\), and therefore
\[
    \frac{D}{Y^{u_2}}
    \geq
    \frac{X^{1-u_2/c_2}}{(\log X)^5}
    \longrightarrow\infty
    \qquad (X\to\infty).
\]
Thus \(D\geq Y^{u_2}\) for all sufficiently large \(X\). Upon choosing
\(c_1\) sufficiently large, all the hypotheses of
Theorem~\ref{thm:fundamental-sieve} are satisfied, and its upper bound
gives
\[
    S(\mathcal A,Y)
    \leq
    5X\prod_{p\leq Y}
    \left(1-\frac{\nu_f(p)}p\right),
\]
which is precisely \eqref{eq:two-linear-form-sieve}.
\end{proof}
\begin{proof}[Proof of Proposition~\ref{prop:good-pair-count}]
It suffices to fix \(\varepsilon\in\{-1,1\}\) and count the solutions to
\begin{equation}\label{eq:signed-linear-relation}
    am-bn=\varepsilon h,
    \qquad
    1\leq h\leq H.
\end{equation}
The two choices of \(\varepsilon\) are treated in the same way, and the
resulting factor of \(2\) is absorbed into the implied constant.

Put
\[
    z:=x^{1/\ell}.
\]
For
\[
    m\in\M_{i,\delta,\ell}(x),
    \qquad
    n\in\M_{j,\delta,\ell}(x),
\]
we write
\[
    m=d_2r_2,\qquad d_2=d_\ell(m),
    \qquad
    n=d_1r_1,\qquad d_1=d_\ell(n).
\]
By the definition of \(\M_{k,\delta,\ell}(x)\),
\[
    d_1,d_2<x^{1/2-\delta}.
\]
Since \(m\) and \(n\) are squarefree, the integers \(d_1,d_2\) are
squarefree and
\[
    (d_1,r_1)=(d_2,r_2)=1.
\]
Moreover, every prime divisor of \(r_1r_2\) is greater than
\(z=x^{1/\ell}\). Since \(r_1,r_2\leq x\), it follows that
\[
    \omega(r_1),\omega(r_2)\leq\ell.
\]
Therefore
\[
    \omega(n)=\omega(d_1)+\omega(r_1),
    \qquad
    \omega(m)=\omega(d_2)+\omega(r_2),
\]
and hence
\[
    i-\ell\leq\omega(d_2)\leq i,
    \qquad
    j-\ell\leq\omega(d_1)\leq j.
\]
For fixed \(d_1,d_2\), equation
\eqref{eq:signed-linear-relation} becomes
\[
    ad_2r_2-bd_1r_1=\varepsilon h.
\]
Set
\[
    g:=(ad_2,bd_1),
    \qquad
    A:=\frac{ad_2}{g},
    \qquad
    B:=\frac{bd_1}{g}.
\]
There are no solutions unless \(g\mid h\). Writing \(h=gq\), we obtain
\begin{equation}\label{eq:normalized-linear-relation}
    Ar_2-Br_1=\varepsilon q,
    \qquad
    (A,B)=1.
\end{equation}
Fix one integral solution
\((r_1^{(0)},r_2^{(0)})\). All integral solutions are then of the form
\begin{equation}\label{eq:parametrized-linear-solutions}
    r_1=r_1^{(0)}+tA,
    \qquad
    r_2=r_2^{(0)}+tB,
    \qquad
    t\in\mathbb Z.
\end{equation}
We also recall
\[
    p\nmid r_1r_2
    \qquad
    \text{for every prime }p\leq z.
\]
For an upper bound, we retain only the restriction \(0<bn\leq x\), thereby
possibly enlarging the set of solutions. Thus every admissible integer
\(t\) lies in the interval
\[
    I
    :=
    \left\{
        t\in\mathbb R:
        0<bd_1\bigl(r_1^{(0)}+tA\bigr)\leq x
    \right\}.
\]
Define
\[
    F(t)
    :=
    \bigl(r_1^{(0)}+tA\bigr)
    \bigl(r_2^{(0)}+tB\bigr).
\]
For every admissible solution, we have
\[
    F(t)=r_1r_2.
\]
Moreover, every prime divisor of \(r_1r_2\) is greater than \(z\).
Consequently, for fixed \(d_1,d_2\) and \(h=gq\), the number of
solutions of \eqref{eq:normalized-linear-relation} satisfying the
original restrictions is at most
\[
    \#\left\{
        t\in I\cap\mathbb Z:
        \left(F(t),\prod_{p\leq z}p\right)=1
    \right\}.
\]
We shall estimate this quantity using
Lemma~\ref{lem:two-linear-form-sieve} with \(X=T\) and \(Y=z\), where
\[
    T:=\#(I\cap\mathbb Z).
\]
Before applying the lemma, we verify the required size hypotheses.

Since \(A>0\), the interval \(I\) has length
\begin{equation}\label{eq:t-interval-length}
    |I|
    =
    \frac{x}{bd_1A}
    =
    \frac{gx}{ab\,d_1d_2}
    =
    \frac{x}{gAB}.
\end{equation}
Using \(g\geq1\), \(ab\leq K^2\), and
\(d_1d_2<x^{1-2\delta}\), we have
\[
    |I|
    \geq
    \frac{x^{2\delta}}{K^2}.
\]
So since
\[
    T=|I|+O(1)
    \qquad\text{and}\qquad
    |I|\geq\frac{x^{2\delta}}{K^2},
\]
we have
\[
    T\geq\frac{x^{2\delta}}{2K^2}
\]
once \(x\) is sufficiently large in terms of \(\delta\) and \(K\).
Since \(\log z=(\log x)/\ell\), it follows that
\[
\begin{aligned}
    \frac{\log T}{\log z}
    &\geq
    \frac{\log\left(x^{2\delta}/(2K^2)\right)}
         {(\log x)/\ell}\\
    &=
    2\delta\ell
    -
    \frac{\ell\log(2K^2)}{\log x}\\
    &=
    2\delta\ell-o_K(1),
\end{aligned}
\]
uniformly for \(\ell\leq L\).

Choose \(\ell_0=\ell_0(\delta)\) sufficiently large that
\[
    2\delta\ell_0>c_2+1,
\]
where \(c_2\) is the absolute constant appearing in
Lemma~\ref{lem:two-linear-form-sieve}. Since \(\ell\leq L\), we also
have
\[
    z=x^{1/\ell}\geq x^{1/L}\ge 2
\] once $x$ is sufficiently large.
It follows that, for every fixed \(K\), uniformly for
\(L\ge \ell\geq\ell_0\),
\[
    T\geq c_1,
    \qquad
    2\leq z\leq T,
    \qquad
    \log T\geq c_2\log z
\]
once \(x\) is sufficiently large.
For each prime \(p\), define
\[
    \nu_F(p)
    :=
    \#\{t\bmod p:F(t)\equiv0\pmod p\}.
\]
Lemma~\ref{lem:two-linear-form-sieve}, applied with \(X=T\) and
\(Y=z\), now gives
\[
\begin{aligned}
    &\#\left\{
        t\in I\cap\mathbb Z:
        \left(F(t),\prod_{p\leq z}p\right)=1
    \right\}\\
    &\qquad\leq
    5T\prod_{p\leq z}
    \left(1-\frac{\nu_F(p)}p\right).
\end{aligned}
\]
We now determine the local factors \(\nu_F(p)\) for primes \(p\leq z\).
Suppose first that \(p\nmid AB\). Both linear factors of \(F\) are
nonconstant modulo \(p\), and their roots coincide if and only if
\(p\mid q\). Hence
\[
    \nu_F(p)
    =
    \begin{cases}
        2,&p\nmid q,\\
        1,&p\mid q.
    \end{cases}
\]
If \(p\mid A\), then \(p\nmid B\). Reducing
\eqref{eq:normalized-linear-relation} modulo \(p\), we find that
\(r_1\) is a nonzero constant modulo \(p\) when \(p\nmid q\), and hence
\(\nu_F(p)=1\). If \(p\mid q\), then \(p\mid r_1\) for every \(t\), so
there are no sifted solutions. The case \(p\mid B\) is symmetric.

Consequently, after discarding configurations whose contribution is
zero, we have
\[
    \nu_F(p)\in\{1,2\},
\]
with \(\nu_F(p)=1\) only if
\[
    p\mid qAB.
\]

The factor corresponding to \(p=2\) is at most \(1\) and may therefore
be omitted. For \(p>2\), replacing \(1-1/p\) by the generic factor \(1-2/p\) introduces the correction
\[
    \frac{1-1/p}{1-2/p}
    =
    1+\frac1{p-2}.
\]
Furthermore,
\[
    \frac{1+1/(p-2)}{1+1/p}
    =
    1+O\left(\frac1{p^2}\right).
\]
It follows that, uniformly over all finite sets of odd primes,
\[
    \prod_p\left(1+\frac1{p-2}\right)
    \asymp
    \prod_p\left(1+\frac1p\right).
\]
Therefore,
\begin{align}\label{eq:sieved-parameter-count}
    &\#\left\{
        t\in I\cap\mathbb Z:
        p\nmid F(t)\text{ for every prime }p\leq z
    \right\} \notag\\
    &\qquad\ll_K
    \frac{gx}{ab\,d_1d_2}
    \prod_{\substack{p\leq z\\p>2}}
    \left(1-\frac2p\right)
    \prod_{p\mid qAB}
    \left(1+\frac1p\right).
\end{align}
For every prime \(p>2\),
\[
    1-\frac2p
    =
    \left(1-\frac1p\right)^2
    \left(1-\frac1{(p-1)^2}\right).
\]
Since
\[
    \prod_{p>2}
    \left(1-\frac1{(p-1)^2}\right)
\]
converges to a strictly positive constant, Mertens' theorem gives
\[
    \prod_{\substack{p\leq z\\p>2}}
    \left(1-\frac2p\right)
    \asymp
    \frac1{(\log z)^2}
    =
    \frac{\ell^2}{(\log x)^2}.
\]

Define
\[
    w(n):=\prod_{p\mid n}\left(1+\frac1p\right).
\]
Since \(w(uv)\leq w(u)w(v)\), summing
\eqref{eq:sieved-parameter-count} over \(h=gq\leq H\) gives
\[
\begin{aligned}
    \sum_{q\leq H/g}w(qAB)
    &\leq
    w(A)w(B)\sum_{q\leq H/g}w(q)\\
    &\ll
    \frac Hg\,w(A)w(B).
\end{aligned}
\]
To see the last inequality, we observe
\[
    w(q)=\sum_{e\mid q}\frac{\mu^2(e)}e
\]
and so, for every \(Q\geq1\),
\[
\begin{aligned}
    \sum_{q\leq Q}w(q)
    &=
    \sum_{e\leq Q}\frac{\mu^2(e)}e
    \left\lfloor\frac Qe\right\rfloor\\
    &\leq
    Q\sum_{e\leq Q}\frac{\mu^2(e)}{e^2}
    \ll Q.
\end{aligned}
\]

Since
\[
    A\mid ad_2,
    \qquad
    B\mid bd_1,
\]
we have
\[
    w(A)\leq w(ad_2)\ll_K w(d_2),
    \qquad
    w(B)\leq w(bd_1)\ll_K w(d_1).
\]
Thus, after summing over \(h\), the contribution associated with fixed
\(d_1,d_2\) is
\[
    \ll_K
    \frac{\ell^2xH}{ab(\log x)^2}
    \frac{w(d_1)w(d_2)}{d_1d_2}.
\]

We now sum over \(d_1\) and \(d_2\). For every \(r\geq0\),
\[
\begin{aligned}
    \sum_{\substack{d<x^{1/2-\delta}\\
                    \mu^2(d)=1\\
                    \omega(d)=r}}
    \frac{w(d)}d
    &\leq
    \frac1{r!}
    \left(
        \sum_{p<x^{1/2-\delta}}
        \frac1p\left(1+\frac1p\right)
    \right)^r\\
    &\ll_{\delta,\gamma}
    \frac{L^r}{r!},
\end{aligned}
\]
uniformly for \(r\leq\gamma L\). Here we used
\[
    \sum_{p<x^{1/2-\delta}}
    \frac1p\left(1+\frac1p\right)
    =
    L+O_\delta(1)
\]
and
\[
    \bigl(L+O_\delta(1)\bigr)^r
    \ll_{\delta,\gamma}
    L^r
    \qquad (r\leq\gamma L).
\]
After dropping all remaining restrictions on \(d_1,d_2\), we obtain
\begin{align}\label{eq:two-slice-factorial-sums}
    &\#\left\{
        (m,n)\in
        \M_{i,\delta,\ell}(x)\times\M_{j,\delta,\ell}(x):
        0<|am-bn|\leq H,\ bn\leq x
    \right\} \notag\\
    &\qquad\ll_{K,\delta,\gamma}
    \frac{\ell^2xH}{ab(\log x)^2}
    \left(
        \sum_{r=i-\ell}^{i}\frac{L^r}{r!}
    \right)
    \left(
        \sum_{s=j-\ell}^{j}\frac{L^s}{s!}
    \right).
\end{align}
 
Because \(L<i<\gamma L\), for \(i-\ell\leq r\leq i\) we have
\[
    \frac{L^r/r!}{L^i/i!}
    =
    \prod_{v=r+1}^{i}\frac{v}{L}
    \leq
    \gamma^{i-r}.
\]
It follows that
\[
    \sum_{r=i-\ell}^{i}\frac{L^r}{r!}
    \ll_\gamma
    \gamma^\ell\frac{L^i}{i!}.
\]
The same argument gives
\[
    \sum_{s=j-\ell}^{j}\frac{L^s}{s!}
    \ll_\gamma
    \gamma^\ell\frac{L^j}{j!}.
\]
 
By Lemma~\ref{lem:sathe-standard}, applied with \(B=\gamma\), we have
uniformly for \(L<i,j<\gamma L\) that
\[
    |N_i(x)|
    \asymp_\gamma
    \frac{x}{\log x}\frac{L^{i-1}}{(i-1)!},
    \qquad
    |N_j(x)|
    \asymp_\gamma
    \frac{x}{\log x}\frac{L^{j-1}}{(j-1)!}.
\]
Moreover, since \(L<i,j<\gamma L\),
\[
    \frac{L^i}{i!}
    =
    \frac{L}{i}\frac{L^{i-1}}{(i-1)!}
    \asymp_\gamma
    \frac{\log x}{x}|N_i(x)|,
\]
and similarly
\[
    \frac{L^j}{j!}
    \asymp_\gamma
    \frac{\log x}{x}|N_j(x)|.
\]
Substituting these estimates into
\eqref{eq:two-slice-factorial-sums}, we obtain
\[
    \ll_{K,\delta,\gamma}
    \frac{\gamma^{2\ell}\ell^2H}{abx}
    |N_i(x)|\,|N_j(x)|.
\]
Taking
\[
    C_\gamma:=\gamma^2
\]
proves \eqref{eq:single-pair-count}.
\end{proof}
 
\section{Two combinatorial tools}\label{sec:prelim}

We record two graph-theoretic and probabilistic tools used to extract a
well-spaced (independent) set of integers from the sparse auxiliary
graphs constructed in Section~\ref{sec:main-proof}. The first, the Caro--Wei
bound, guarantees an independent set of size at least
$\sum_{v\in V(G)}(d_G(v)+1)^{-1}$; combined with the
Cauchy--Schwarz inequality, it converts
an upper bound on the number of edges of a graph into a lower bound on the
size of an independent set, and is used in parts~\textup{(i)} and
\textup{(iii)} of the main theorem. Part~\textup{(ii)} requires more: an
independent set that is simultaneously large on \emph{every} slice
\(\{n:\omega(n)=k\}\) with \(k\) in a range of width
\(\asymp\sqrt L\), rather than merely large in total. For this we use a
bounded-differences inequality for functions on the symmetric group under
arbitrary transpositions. It lets us choose a single permutation realizing large independent sets on all the slices at once, via a union bound.

\begin{theorem}[Caro--Wei {\cite{AlonSpencer}}]\label{cw}
    Given an undirected, finite graph \(G\) with vertex set \(V(G)\), there is
an independent set \(I\) of \(G\) with
    \begin{equation*}
        |I|\ge \sum_{v \in V(G)}\frac{1}{d_G(v)+1}.
    \end{equation*}
\end{theorem}

\begin{theorem}[Permutation bounded-differences inequality {\cite[Theorem~3.1]{DudekGrytczukRucinski}}]\label{ah}
    Let $n \in \mathbb{N}$ and let $h: \mathbb{S}_n\rightarrow \mathbb{R}$ be such that for every $\sigma \in \mathbb{S}_n$ and every transposition $\pi\in \mathbb{S}_n$,
    \begin{equation*}
        |h(\sigma)-h(\sigma\circ\pi)|\le c
    \end{equation*}
    for some constant $c>0$. If $X$ is a uniformly random permutation in $\mathbb{S}_n$, then for every \(t>0,\)
    \begin{equation*}
        \mathbb{P}\bigl(|h(X)-\mathbb{E}[h(X)]|>t\bigr)
        \le 2\exp\left(-\dfrac{t^2}{2nc^2}\right).
    \end{equation*}
\end{theorem}

\section{Proof of the main theorem}\label{sec:main-proof}

We now prove Theorem~\ref{thm:main}. 
Recall that
\[
    L:=\log\log x,
    \qquad
    I(t):=1+t(\log t-1),
    \qquad
    y_t:=(\log x)^{I(t)}.
\]

\subsection{The supercritical regime
\texorpdfstring{\(\lambda>2\)}{lambda > 2}}
\label{subsec:supercritical}

\begin{proof}[Proof of Theorem~\ref{thm:main}\textup{(i)}]
Fix \(\lambda>2\), and put
\[
    y:=y_\lambda=(\log x)^{I(\lambda)},
    \qquad
    T_0:=\frac{\log L}{2\log\lambda}.
\]
We prove the upper and lower bounds for \(S(x,y)\) in turn.

\smallskip
\emph{Upper bound.}
Let \(\mathcal S\subseteq[1,x]\) be such that
\[
    |a-b|>y
    \qquad(a,b\in\mathcal S,\ a\neq b).
\]
Ordering the elements of \(\mathcal S\) increasingly, consecutive
elements differ by more than \(y\). Hence
\begin{equation}\label{eq:spaced-cardinality}
    |\mathcal S|<\frac{x}{y}+1.
\end{equation}
We split
\begin{equation}\label{eq:supercritical-split}
\begin{aligned}
    \sum_{n\in\mathcal S}\tau(n)
    &=
    \sum_{\substack{n\in\mathcal S\\
          \tau(n)>2^{\lambda L-T_0}}}\tau(n)
    +
    \sum_{\substack{n\in\mathcal S\\
          \tau(n)\leq 2^{\lambda L-T_0}}}\tau(n)\\
    &=:\mathcal S_1+\mathcal S_2.
\end{aligned}
\end{equation}
We first bound \(\mathcal S_1\).

We recall that every positive integer \(n\) has a unique decomposition
\[
    n=m(n)r(n),
\]
where \(m(n)\) is powerful, \(r(n)\) is squarefree, and
\((m(n),r(n))=1\). Consequently,
\[
    \tau(n)=\tau(m(n))\tau(r(n))
           =\tau(m(n))2^{\omega(r(n))}.
\]
It follows that
\begin{equation}\label{eq:supercritical-high-tail}
\begin{aligned}
    \mathcal S_1
    &\leq
    \sum_{\substack{n\leq x\\
          \tau(n)>2^{\lambda L-T_0}}}\tau(n)\\
    &\leq
    \sum_{\substack{m\leq x\\m\ \mathrm{powerful}}}\tau(m)
    \left(
    \sum_{\substack{d\leq x/m\\
          \omega(d)>
          \lambda L-T_0-\frac{\log\tau(m)}{\log 2}}}
          \mu^2(d)\tau(d)
    \right).
\end{aligned}
\end{equation}
The coprimality condition \((m,d)=1\), which is present in the
canonical decomposition, has been discarded in the
last upper bound.

In \eqref{eq:supercritical-high-tail}, split the sum over powerful
\(m\) at \(m=x^{1/2}\):
\begin{equation}\label{eq:supercritical-powerful-split}
\begin{aligned}
    &\sum_{\substack{m\leq x\\m\ \mathrm{powerful}}}\tau(m)
    \left(
    \sum_{\substack{d\leq x/m\\
          \omega(d)>
          \lambda L-T_0-\frac{\log\tau(m)}{\log 2}}}
          \mu^2(d)\tau(d)
    \right)\\
    &\quad=
    \sum_{\substack{m\leq x^{1/2}\\m\ \mathrm{powerful}}}
       \tau(m)
    \left(
    \sum_{\substack{d\leq x/m\\
          \omega(d)>
          \lambda L-T_0-\frac{\log\tau(m)}{\log 2}}}
          \mu^2(d)\tau(d)
    \right)\\
    &\qquad+
    \sum_{\substack{x^{1/2}<m\leq x\\m\ \mathrm{powerful}}}
       \tau(m)
    \left(
    \sum_{\substack{d\leq x/m\\
          \omega(d)>
          \lambda L-T_0-\frac{\log\tau(m)}{\log 2}}}
          \mu^2(d)\tau(d)
    \right)\\
    &\quad=: \Sigma_1+\Sigma_2.
\end{aligned}
\end{equation}

We deal with \(\Sigma_2\) first. Recall that every powerful integer
has a representation
\[
    m=a^2b^3,
\]
where \(b\) is squarefree. To see this, write
\[
    m=\prod_p p^{e_p},
    \qquad e_p\geq2.
\]
Every even exponent is absorbed into \(a^2\), whereas every odd
exponent \(e_p\geq3\) is written as
\[
    e_p=2\left(\frac{e_p-3}{2}\right)+3,
\]
with the corresponding prime included in \(b\). Thus, if
\[
    P(U):=\#\{m\leq U:m\ \mathrm{powerful}\},
\]
then, on extending the sum from squarefree \(b\) to all positive
integers \(b\),
\begin{equation}\label{eq:powerful-count}
\begin{aligned}
    P(U)
    &\leq
    \sum_{b\leq U^{1/3}}
       \left(\frac{U}{b^3}\right)^{1/2}\\
    &\ll
    U^{1/2}\sum_{b\geq1}\frac1{b^{3/2}}
    \ll U^{1/2}.
\end{aligned}
\end{equation}
We also have, for every \(\eta>0\),
\[
    \tau(m)\ll_\eta m^\eta.
\]
Therefore, by partial summation, for every fixed
\(0<\eta<1/2\),
\begin{align}
    \sum_{\substack{m>U\\m\ \mathrm{powerful}}}
       \frac{\tau(m)}m
    &\ll_\eta
    \sum_{\substack{m>U\\m\ \mathrm{powerful}}}
       m^{-1+\eta}                                                     \notag\\
    &\ll_\eta
    \frac{P(U)}{U^{1-\eta}}
    +\int_U^\infty\frac{P(t)}{t^{2-\eta}}\,dt                         \notag\\
    &\ll_\eta U^{-1/2+\eta}.
    \label{eq:powerful-tail}
\end{align}
The integral converges because \(\eta<1/2\). Taking
\(U=x^{1/2}\) in \eqref{eq:powerful-tail} gives
\begin{equation}\label{eq:powerful-tail-sqrtx}
    \sum_{\substack{m>x^{1/2}\\m\ \mathrm{powerful}}}
       \frac{\tau(m)}m
    \ll_\eta x^{-1/4+\eta/2}.
\end{equation}
For the inner sum in \(\Sigma_2\), we use the unrestricted estimate
\[
    \sum_{d\leq T}\mu^2(d)\tau(d)
    \leq\sum_{d\leq T}\tau(d)
    \ll T\log(2T).
\]
Consequently,
\begin{equation}\label{eq:supercritical-large-powerful}
\begin{aligned}
    \Sigma_2
    &\ll
    x\log x
    \sum_{\substack{m>x^{1/2}\\m\ \mathrm{powerful}}}
       \frac{\tau(m)}m\\
    &\ll_\eta x^{3/4+\eta/2}\log x.
\end{aligned}
\end{equation}
Since \(\lambda>2\) is fixed and
\[
    \frac{x(\log x)^{\lambda\log 2}}
         {yL^{\log 2/(2\log\lambda)}}
\]
is \(x\) multiplied by fixed powers of \(\log x\) and \(L\),
\eqref{eq:supercritical-large-powerful} and \(\eta<1/2\) imply
\begin{equation}\label{eq:supercritical-large-powerful-negligible}
    \Sigma_2
    =o_\lambda\left(
       \frac{x(\log x)^{\lambda\log 2}}
            {yL^{\log 2/(2\log\lambda)}}
    \right).
\end{equation}
Indeed, after division by the quantity in parentheses, the bound in
\eqref{eq:supercritical-large-powerful} contains the factor
\(x^{-1/4+\eta/2}\), which dominates every fixed power of
\(\log x\) and \(L\).

We now estimate \(\Sigma_1\). Since \(\lambda>2\) is fixed, put
\[
    \delta_0:=\frac{\lambda-2}{2}.
\]
Then there is a constant \(c_\lambda>0\), depending only on
\(\lambda\), such that, for all sufficiently large \(x\),
\begin{equation}\label{eq:supercritical-delta0}
\begin{aligned}
    \lambda L-T_0-\delta_0L
    &=(\lambda-\delta_0)L-T_0\\
    &=(\lambda-\delta_0)L
      -\frac{\log L}{2\log\lambda}\\
    &>(2+c_\lambda)L.
\end{aligned}
\end{equation}
We split \(\Sigma_1\) according as
\(\tau(m)\leq2^{\delta_0L}\) or
\(\tau(m)>2^{\delta_0L}\):
\begin{equation}\label{eq:supercritical-sigma1-split}
\begin{aligned}
    \Sigma_1
    &=
    \sum_{\substack{m\leq x^{1/2}\\
          m\ \mathrm{powerful}\\
          \tau(m)\leq2^{\delta_0L}}}
       \tau(m)
    \left(
    \sum_{\substack{d\leq x/m\\
          \omega(d)>
          \lambda L-T_0-\frac{\log\tau(m)}{\log 2}}}
          \mu^2(d)\tau(d)
    \right)\\
    &\qquad+
    \sum_{\substack{m\leq x^{1/2}\\
          m\ \mathrm{powerful}\\
          \tau(m)>2^{\delta_0L}}}
       \tau(m)
    \left(
    \sum_{\substack{d\leq x/m\\
          \omega(d)>
          \lambda L-T_0-\frac{\log\tau(m)}{\log 2}}}
          \mu^2(d)\tau(d)
    \right)\\
    &=:\Gamma_1+\Gamma_2.
\end{aligned}
\end{equation}

We first treat \(\Gamma_1\). For a fixed powerful
\(m\leq x^{1/2}\), put
\[
    T:=\frac{x}{m},
    \qquad
    L_T:=\log\log T.
\]
Since \(x^{1/2}\leq T\leq x\), we have
\[
    L-\log 2\leq L_T\leq L;
\]
in particular, \(L_T=L+O(1)\), uniformly for
\(m\leq x^{1/2}\).

The summation over \(\omega(d)\) extends through the entire upper
tail. We therefore use the all-\(j\) Hardy--Ramanujan bound from
Lemma~\ref{lem:hardy-ramanujan}, rather than the bounded-range
Sathe--Selberg asymptotic in Lemma~\ref{lem:sathe-standard}. Let
\(c_2\) be the absolute constant in the Hardy--Ramanujan bound, and
put
\[
    A_m:=\lambda L-T_0-\frac{\log\tau(m)}{\log 2},
    \qquad
    R_m:=\lfloor A_m\rfloor.
\]
Since \(\tau(m)\leq2^{\delta_0L}\), we may decrease
\(c_\lambda\), if necessary, and use
\eqref{eq:supercritical-delta0} to obtain
\[
    R_m>(2+c_\lambda/2)L
\]
for all sufficiently large \(x\). The Hardy--Ramanujan bound gives
\begin{align}
    &\sum_{\substack{d\leq T\\
          \omega(d)>A_m}}
          \mu^2(d)\tau(d)                                      \notag\\
    &\quad=
    \sum_{j>A_m}2^j|N_j(T)|                                    \notag\\
    &\quad\ll
    \frac{T}{\log T}
    \sum_{j>A_m}
       2^j\frac{(L_T+c_2)^{j-1}}{(j-1)!}                       \notag\\
    &\quad=
    \frac{2T}{\log T}
    \sum_{h\geq R_m}
       \frac{\{2(L_T+c_2)\}^{h}}{h!}.
    \label{eq:supercritical-HR-expansion}
\end{align}
The equality in the last line follows by putting \(h=j-1\): the
least integer \(j>A_m\) is \(R_m+1\), so the least corresponding
value of \(h\) is \(R_m\).

For \(h\geq R_m\), the ratio of two consecutive summands in the
last factorial series is
\[
    \frac{2(L_T+c_2)}{h+1}.
\]
Since \(R_m>(2+c_\lambda/2)L\) and \(L_T=L+O(1)\), this ratio is
at most a fixed number \(\rho_\lambda<1\), once \(x\) is large.
The series is therefore bounded by a constant multiple of its first
term. Moreover, \(R_m=O_\lambda(L)\), and hence
\[
    \left(\frac{L_T+c_2}{L}\right)^{R_m}=O_\lambda(1).
\]
Thus replacing \(L_T+c_2\) by \(L\) costs only a bounded
\(\lambda\)-dependent factor, and
\begin{equation}\label{eq:HR-tail-first-term}
    \sum_{\substack{d\leq T\\
          \omega(d)>A_m}}
          \mu^2(d)\tau(d)
    \ll_\lambda
    \frac{x}{m\log x}\frac{(2L)^{R_m}}{R_m!}.
\end{equation}
Here we also used \(T=x/m\) and
\(\log T\asymp\log x\), uniformly for \(m\leq x^{1/2}\).

Let
\[
    K_0:=\lfloor\lambda L\rfloor,
    \qquad
    \beta_\lambda:=\frac{\log(\lambda/2)}{\log 2}.
\]
Then
\[
    K_0-R_m
    =T_0+\frac{\log\tau(m)}{\log 2}+O(1).
\]
Since \(R_m\leq K_0\), we have
\begin{align}
    \frac{(2L)^{R_m}}{R_m!}
    &=
    \frac{(2L)^{K_0}}{K_0!}
    \prod_{\nu=R_m+1}^{K_0}\frac{\nu}{2L}                       \notag\\
    &\leq
    \frac{(2L)^{K_0}}{K_0!}
    \left(\frac{\lambda}{2}\right)^{K_0-R_m}                   \notag\\
    &\ll_\lambda
    \frac{(2L)^{K_0}}{K_0!}
    \left(\frac{\lambda}{2}\right)^{T_0}
    \tau(m)^{\beta_\lambda}.
    \label{eq:supercritical-factorial-comparison}
\end{align}
The exponent \(\beta_\lambda\) is fixed. We may use the cruder
inequality
\[
    \beta_\lambda\leq\frac{\lambda\log 2}{2};
\]
indeed, with \(u=\lambda/2\), this follows from
\(\sup_{u>0}(\log u)/u=1/e<(\log 2)^2\). Therefore
\begin{equation}\label{eq:supercritical-factorial-comparison-crude}
    \frac{(2L)^{R_m}}{R_m!}
    \ll_\lambda
    \frac{(2L)^{K_0}}{K_0!}
    \left(\frac{\lambda}{2}\right)^{T_0}
    \tau(m)^{\lambda\log 2/2}.
\end{equation}

Now invoking the inequalities \eqref{eq:robbins-stirling} and recalling that \(K_0=\lambda L+O(1)\), these inequalities give
\begin{align*}
    \frac{x}{\log x}\frac{(2L)^{K_0}}{K_0!}
    &\asymp_\lambda
    \frac{x}{\log x\sqrt L}
       \left(\frac{2e}{\lambda}\right)^{\lambda L}\\
    &=
    \frac{x}{\sqrt L}
       (\log x)^{-1+\lambda(1+\log 2-\log\lambda)}\\
    &=
    \frac{x(\log x)^{\lambda\log 2}}
         {y\sqrt L},
\end{align*}
where in the last line we used
\(I(\lambda)=1-\lambda+\lambda\log\lambda\). Also,
\[
    \left(\frac{\lambda}{2}\right)^{T_0}
    =L^{(\log\lambda-\log 2)/(2\log\lambda)}
    =L^{1/2-\log 2/(2\log\lambda)}.
\]
Combining these estimates with
\eqref{eq:HR-tail-first-term} and
\eqref{eq:supercritical-factorial-comparison-crude}, we obtain
\begin{equation}\label{eq:supercritical-gamma1}
\begin{aligned}
    \Gamma_1
    &\ll_\lambda
    \frac{x(\log x)^{\lambda\log 2}}
         {yL^{\log 2/(2\log\lambda)}}
    \sum_{\substack{m\leq x^{1/2}\\m\ \mathrm{powerful}}}
       \frac{\tau(m)^{1+\lambda\log 2/2}}m\\
    &\ll_\lambda
    \frac{x(\log x)^{\lambda\log 2}}
         {yL^{\log 2/(2\log\lambda)}}.
\end{aligned}
\end{equation}
For completeness, the convergence used in the last step follows from
the fact that, for every fixed \(q>0\),
\begin{equation}\label{eq:powerful-moment}
\begin{aligned}
    \sum_{\substack{m\geq1\\m\ \mathrm{powerful}}}
       \frac{\tau(m)^q}{m}
    &=
    \prod_p\left(
       1+\sum_{\nu\geq2}\frac{(\nu+1)^q}{p^\nu}
    \right)\\
    &<\infty,
\end{aligned}
\end{equation}
because every local factor is \(1+O_q(p^{-2})\).

It remains to bound \(\Gamma_2\). For any \(U>1\) and any fixed
\(q>1\), Rankin's inequality and \eqref{eq:powerful-moment} give
\begin{equation}\label{eq:supercritical-rankin}
\begin{aligned}
    \sum_{\substack{m\geq1\\m\ \mathrm{powerful}\\
                    \tau(m)>U}}
       \frac{\tau(m)}m
    &\leq
    U^{-(q-1)}
    \sum_{\substack{m\geq1\\m\ \mathrm{powerful}\\
                    \tau(m)>U}}
       \frac{\tau(m)^q}{m}\\
    &\ll_q U^{-(q-1)}.
\end{aligned}
\end{equation}
Apply this with \(U=2^{\delta_0L}\), and use once more the
unrestricted divisor-sum estimate for the inner sum in
\(\Gamma_2\). We obtain
\begin{align}
    \Gamma_2
    &\ll
    x\log x
    \sum_{\substack{m\leq x^{1/2}\\m\ \mathrm{powerful}\\
                    \tau(m)>2^{\delta_0L}}}
       \frac{\tau(m)}m                                           \notag\\
    &\ll_q
    x\log x\,2^{-\delta_0(q-1)L}.                                \notag
\end{align}
Given \(A>0\), choose the fixed number
\(q=q(A,\delta_0)>1\) so that
\[
    \delta_0(q-1)\log 2>A+1.
\]
Since
\[
    2^{\delta_0L}=(\log x)^{\delta_0\log 2},
\]
it follows that
\begin{equation}\label{eq:supercritical-gamma2}
    \Gamma_2\ll_{A,\lambda}\frac{x}{(\log x)^A}
\end{equation}
for every fixed \(A>0\). Finally, because
\[
    y=(\log x)^{I(\lambda)}
     =(\log x)^{1-\lambda+\lambda\log\lambda},
\]
we may choose \(A\), depending only on \(\lambda\), sufficiently
large that
\[
    A>I(\lambda)-\lambda\log 2.
\]
For this choice,
\[
    \frac{x}{(\log x)^A}
    =o_\lambda\left(
       \frac{x(\log x)^{\lambda\log 2}}
            {yL^{\log 2/(2\log\lambda)}}
    \right).
\]
Thus
\begin{equation}\label{eq:supercritical-large-tau-powerful}
    \Gamma_2
    \ll_\lambda
    \frac{x(\log x)^{\lambda\log 2}}
         {yL^{\log 2/(2\log\lambda)}}.
\end{equation}

Combining \eqref{eq:supercritical-gamma1} and
\eqref{eq:supercritical-large-tau-powerful}, we have
\begin{equation}\label{eq:supercritical-sigma1}
    \Sigma_1
    \ll_\lambda
    \frac{x(\log x)^{\lambda\log 2}}
         {yL^{\log 2/(2\log\lambda)}}.
\end{equation}
Hence \eqref{eq:supercritical-high-tail},
\eqref{eq:supercritical-powerful-split}, and
\eqref{eq:supercritical-large-powerful-negligible} yield
\begin{equation}\label{eq:supercritical-unrestricted-high-tail}
    \sum_{\substack{n\leq x\\
          \tau(n)>2^{\lambda L-T_0}}}\tau(n)
    \ll_\lambda
    \frac{x(\log x)^{\lambda\log 2}}
         {yL^{\log 2/(2\log\lambda)}}.
\end{equation}
In particular, this bounds \(\mathcal S_1\). Notice that
\eqref{eq:supercritical-unrestricted-high-tail} concerns only the
high-divisor tail; it is not an estimate for the unrestricted sum
\(\sum_{n\leq x}\tau(n)\).

For the low-divisor portion, \eqref{eq:spaced-cardinality} gives
\begin{align*}
    \mathcal S_2
    &\leq
    2^{\lambda L-T_0}|\mathcal S|\\
    &<
    2^{\lambda L-T_0}\left(\frac{x}{y}+1\right).
\end{align*}
Now
\[
    2^{\lambda L-T_0}
    =\frac{(\log x)^{\lambda\log 2}}
           {L^{\log 2/(2\log\lambda)}}.
\]
The term containing \(x/y\) is therefore exactly of the desired
order, while the additional term is negligible because
\[
    \frac{2^{\lambda L-T_0}}
         {x(\log x)^{\lambda\log 2}/
          \{yL^{\log 2/(2\log\lambda)}\}}
    =\frac{y}{x}=o_\lambda(1).
\]
It follows that
\begin{equation}\label{eq:supercritical-S2}
    \mathcal S_2
    \ll_\lambda
    \frac{x(\log x)^{\lambda\log 2}}
         {yL^{\log 2/(2\log\lambda)}}.
\end{equation}
Together with \eqref{eq:supercritical-unrestricted-high-tail}, this
gives
\[
    \sum_{n\in\mathcal S}\tau(n)
    \ll_\lambda
    \frac{x(\log x)^{\lambda\log 2}}
         {yL^{\log 2/(2\log\lambda)}}.
\]
Since \(\mathcal S\) was an arbitrary subset of \([1,x]\) whose
distinct elements have mutual distance greater than \(y\), taking the
supremum over all such \(\mathcal S\) gives
\begin{equation}\label{eq:supercritical-upper}
    S(x,y)
    \ll_\lambda
    \frac{x(\log x)^{\lambda\log 2}}
         {y(\log\log x)^{\log 2/(2\log\lambda)}}.
\end{equation}
This proves the upper bound.

\smallskip
\emph{Lower bound.}
Fix \(\delta=1/4\), and put
\[
    k:=\lfloor\lambda L-T_0\rfloor.
\]
For all sufficiently large \(x\), we have \(L<k<\lambda L\).
Choose a fixed integer \(\ell\) large enough that
\(\ell\geq\ell_0(\delta)\), as required in
Proposition~\ref{prop:good-pair-count}, and large enough in terms of
\(\lambda\) and \(\delta\) for the final assertion of
Proposition~\ref{prop:typical-useful-divisor}\textup{(a)} to apply
with \(\gamma=\lambda\). Since \(\ell\) is fixed, \(\ell<k\) and
\(\ell\leq L\) also hold for all sufficiently large \(x\).
The useful-divisor proposition gives
\begin{equation}\label{eq:supercritical-typical-set}
    |\M_{k,\delta,\ell}(x)|
    \gg_{\lambda,\delta}|N_k(x)|.
\end{equation}
We also set
\[
    K_0:=\lfloor\lambda L\rfloor.
\]
We first record the size of this Sathe--Selberg layer, since the shift
by \(T_0\) is chosen precisely so that its cardinality is
\(\asymp x/y\). By Robbins' explicit factorial inequalities
\eqref{eq:robbins-stirling},
\begin{equation}\label{eq:supercritical-central-layer}
    \frac{x}{\log x}\frac{L^{K_0}}{K_0!}
    \asymp_\lambda\frac{x}{y\sqrt L}.
\end{equation}
Indeed, the same calculation used in the upper bound gives
\[
    \frac{x}{\log x}\frac{L^{K_0}}{K_0!}
    \asymp_\lambda
    \frac{x}{\log x\sqrt L}
       \left(\frac{e}{\lambda}\right)^{\lambda L}
    =\frac{x}{y\sqrt L}.
\]
Moreover,
\[
    K_0-k=T_0+O(1),
\]
and hence
\begin{align*}
    \frac{L^k/k!}{L^{K_0}/K_0!}
    &=L^{-(K_0-k)}\frac{K_0!}{k!}\\
    &=\prod_{\nu=k+1}^{K_0}\frac{\nu}{L}\\
    &=\lambda^{K_0-k}
      \exp\left(
         O_\lambda\left(\frac{(K_0-k)^2}{L}\right)
      \right)\\
    &\asymp_\lambda\lambda^{T_0}
     =\sqrt L.
\end{align*}
Here \((K_0-k)^2/L=O((\log L)^2/L)=o(1)\), and
\[
    \lambda^{T_0}
    =\exp(T_0\log\lambda)=\sqrt L.
\]
Lemma~\ref{lem:sathe-standard}, applied with $B=\lambda$ now gives
\[
    |N_k(x)|
    \asymp_\lambda
    \frac{x}{\log x}\frac{L^{k-1}}{(k-1)!}.
\]
The shift between \(L^{k-1}/(k-1)!\) and \(L^k/k!\) is harmless,
because their ratio is \(k/L\asymp_\lambda 1\). Therefore
\begin{equation}\label{eq:supercritical-Nk}
    |N_k(x)|\asymp_\lambda\frac{x}{y}.
\end{equation}

We next apply Proposition~\ref{prop:good-pair-count} with $\gamma=\lambda+1$
and take
\[
    K=1,
    \qquad a=b=1,
    \qquad i=j=k,
    \qquad H=y.
\]
For all sufficiently large \(x\), we have
\(L<k<\gamma L\), \(\ell\leq L\), and \(y\geq1\), so all the
hypotheses of the proposition are satisfied. Hence
\begin{equation}\label{eq:supercritical-pair-count}
\begin{aligned}
    &\#\bigl\{(m,n)\in
       \M_{k,\delta,\ell}(x)^2:
       0<|m-n|\leq y,\ n\leq x\bigr\}\\
    &\qquad\ll_{\lambda,\delta}
       C_{\lambda+1}^\ell\ell^2\frac{y}{x}|N_k(x)|^2\\
    &\qquad\ll_{\lambda,\delta,\ell}\frac{x}{y},
\end{aligned}
\end{equation}
where the last estimate follows from
\eqref{eq:supercritical-Nk}.

Now form a graph \(G\) with vertex set
\[
    V(G):=\M_{k,\delta,\ell}(x),
\]
joining two distinct vertices \(g\) and \(h\) precisely when
\[
    0<|g-h|\leq y.
\]
Let \(E\) denote the number of its edges. Because \(a=b=1\), every
edge \(\{g,h\}\), with \(g<h\), contributes the
ordered pair \((h,g)\) to the set counted on the left-hand side of
\eqref{eq:supercritical-pair-count}; its second coordinate is
\(g<x\). Consequently,
\begin{equation}\label{eq:supercritical-edge-count}
    E
    \leq
    \#\bigl\{(m,n)\in
       \M_{k,\delta,\ell}(x)^2:
       0<|m-n|\leq y,\ n\leq x\bigr\}
    \ll_{\lambda,\delta,\ell}\frac{x}{y}.
\end{equation}

By the Caro--Wei bound, Theorem~\ref{cw}, the graph \(G\) contains
an independent set \(\mathcal T\) with
\begin{align}
    |\mathcal T|
    &\geq
    \sum_{v\in V(G)}\frac1{d_G(v)+1}                              \notag\\
    &\geq
    \frac{|V(G)|^2}{|V(G)|+2E}                                    \notag\\
    &=
    \frac{|\M_{k,\delta,\ell}(x)|^2}
         {|\M_{k,\delta,\ell}(x)|+2E}                             \notag\\
    &\gg_{\lambda,\delta,\ell}|N_k(x)|
     \gg_\lambda\frac{x}{y}.
    \label{eq:supercritical-independent-set}
\end{align}
The second inequality is Cauchy--Schwarz, and the penultimate one
uses \eqref{eq:supercritical-typical-set},
\eqref{eq:supercritical-Nk}, and
\eqref{eq:supercritical-edge-count}.

Every \(n\in\mathcal T\) is squarefree and satisfies
\(\omega(n)=k\); hence \(\tau(n)=2^k\). It follows that
\begin{align}
    S(x,y)
    &\geq\sum_{n\in\mathcal T}\tau(n)                              \notag\\
    &=2^k|\mathcal T|                                               \notag\\
    &\gg_\lambda 2^{\lambda L-T_0}\frac{x}{y}                     \notag\\
    &\asymp_\lambda
    \frac{x(\log x)^{\lambda\log 2}}
         {y(\log\log x)^{\log 2/(2\log\lambda)}}.
    \label{eq:supercritical-lower}
\end{align}
Finally, since \(\mathcal T\) is independent in \(G\),
\[
    |a-b|>y
    \qquad(a,b\in\mathcal T,\ a\neq b),
\]
so \(\mathcal T\) is admissible in the definition of \(S(x,y)\).
Combining \eqref{eq:supercritical-upper} and
\eqref{eq:supercritical-lower} proves
\[
    S(x,y)
    \asymp_\lambda
    \frac{x(\log x)^{\lambda\log 2}}
         {y(\log\log x)^{\log 2/(2\log\lambda)}},
\]
as required.
\end{proof}


\subsection{The critical transition window
\texorpdfstring{\(1<\lambda\leq2+2\kappa/\sqrt L\)}
{1 < lambda <= 2 + 2 kappa/sqrt(L)}}
\label{subsec:transition}

\begin{proof}[Proof of Theorem~\ref{thm:main}\textup{(ii)}]
Fix $\kappa>0$, and put
\[
    \lambda_*:=2+\frac{2\kappa}{\sqrt L},
    \qquad
    y_*:=y_{\lambda_*},
    \qquad
    T_*:=\frac{\log L}{2\log\lambda_*}.
\]
Since $I'(t)=\log t>0$ for $t>1$, the function
$t\mapsto y_t$ is increasing. Since $S(x,y)$ is non-increasing
in $y$, it suffices to prove
\begin{equation}\label{eq:transition-endpoint}
    S(x,y_*)\gg_\kappa x\log x.
\end{equation}

Define
\[
    \mathscr K
    :=
    \left\{
       k\in\mathbb N:
       \lambda_*L-T_*<k<
       \lambda_*L-T_*+\kappa\sqrt L
    \right\},
\]
and let
\[
    k_0:=\min\mathscr K.
\]
Then
\begin{equation}\label{eq:transition-k0}
    k_0=\lambda_*L-T_*+O(1)
\end{equation}
and
\begin{equation}\label{eq:transition-number-layers}
    |\mathscr K|=\kappa\sqrt L+O(1).
\end{equation}
Since $T_*=O(\log L)$, for all sufficiently large $x$,
\begin{equation}\label{eq:transition-layer-range}
    2L<k<3L
    \qquad(k\in\mathscr K).
\end{equation}

For $k\geq1$, write
\[
    U_k:=\frac{x}{\log x}\frac{L^{k-1}}{(k-1)!}.
\]
By Lemma~\ref{lem:sathe-standard} for $B=3,$ uniformly for
$k\in\mathscr K$,
\begin{equation}\label{eq:transition-sathe}
    |N_k(x)|\asymp U_k.
\end{equation}

Let
\[
    q_*:=\lfloor\lambda_*L\rfloor.
\]
Applying Robbins' inequalities
\eqref{eq:robbins-stirling} to $(q_*-1)!$ gives
\begin{equation}\label{eq:transition-sathe-central}
    U_{q_*}
    \asymp_\kappa
    \frac{x}{\sqrt L}e^{-I(\lambda_*)L}
    =
    \frac{x}{y_*\sqrt L}.
\end{equation}
Moreover,
\[
    q_*-k_0=T_*+O(1),
\]
and hence
\begin{align*}
    \frac{U_{k_0}}{U_{q_*}}
    &=
    \prod_{r=k_0}^{q_*-1}\frac rL\\
    &=
    \exp\left(
       (T_*+O(1))\log\lambda_*
       +O_\kappa\left(\frac{T_*^2}{L}\right)
    \right)\\
    &\asymp_\kappa\lambda_*^{T_*}
     =\sqrt L.
\end{align*}
Thus
\begin{equation}\label{eq:transition-first-layer}
    U_{k_0}\asymp_\kappa\frac{x}{y_*},
    \qquad
    |N_{k_0}(x)|\asymp_\kappa\frac{x}{y_*}.
\end{equation}

For $j\geq k_0$, \eqref{eq:transition-layer-range} gives
\[
    \frac{U_{j+1}}{U_j}=\frac Lj<\frac12.
\]
Therefore
\begin{equation}\label{eq:transition-unweighted-mass}
    \sum_{j\in\mathscr K}|N_j(x)|
    \ll_\kappa\frac{x}{y_*}.
\end{equation}

If $k\in\mathscr K$, then
\[
    0\leq k-k_0\leq\kappa\sqrt L+O(1).
\]
Since every intervening index is less than $3L$,
\[
    \frac{U_k}{U_{k_0}}
    =
    \prod_{r=k_0}^{k-1}\frac Lr
    \geq3^{-(k-k_0)}.
\]
Consequently,
\begin{equation}\label{eq:transition-minimum-layer}
    |N_k(x)|
    \gg_\kappa
    \frac{x}{y_*}\,3^{-\kappa\sqrt L+O(1)}
    \qquad(k\in\mathscr K).
\end{equation}

Fix $\delta\in(0,1/2)$ and choose an integer $\ell$, independently
of $x$, which is sufficiently large that
Proposition~\ref{prop:typical-useful-divisor}\textup{(a)},
applied with $\gamma=3$, gives
\[
    |\M_{k,\delta,\ell}(x)|
    \geq c_0|N_k(x)|
    \qquad(k\in\mathscr K)
\]
for some fixed $c_0>0$. We also choose
$\ell\geq\ell_0(\delta)$ so that
Proposition~\ref{prop:good-pair-count} is applicable.

Let $G$ be the graph with vertex set
\[
    V:=\bigcup_{k\in\mathscr K}\M_{k,\delta,\ell}(x),
\]
where two distinct vertices are adjacent if their distance is at
most $y_*$. For every $k\in\mathscr K$, Proposition~\ref{prop:good-pair-count},
applied with $\gamma=3$ and with
\[
    K=1,\qquad a=b=1,\qquad H=y_*,
\]
gives
\begin{align*}
    \sum_{m\in\M_{k,\delta,\ell}(x)}d_G(m)
    &=
    \sum_{j\in\mathscr K}
    \#\left\{
       (m,n)\in
       \M_{k,\delta,\ell}(x)\times
       \M_{j,\delta,\ell}(x):
       0<|m-n|\leq y_*
    \right\}\\
    &\ll_\delta
    C_3^\ell\ell^2\frac{y_*}{x}
    |N_k(x)|
    \sum_{j\in\mathscr K}|N_j(x)|\\
    &\ll_{\kappa,\delta,\ell}|N_k(x)|.
\end{align*}
Thus there is a constant $B=B(\kappa,\delta,\ell)$ such that
\begin{equation}\label{eq:transition-layer-degree}
    \sum_{m\in\M_{k,\delta,\ell}(x)}d_G(m)
    \leq B|N_k(x)|
    \qquad(k\in\mathscr K).
\end{equation}

Set
\[
    D:=\left\lceil\frac{2B}{c_0}\right\rceil
\]
and
\[
    \M'_k
    :=
    \{m\in\M_{k,\delta,\ell}(x):d_G(m)\leq D\}.
\]
By Markov's inequality on \eqref{eq:transition-layer-degree},
\[
    |\M_{k,\delta,\ell}(x)\setminus\M'_k|
    \leq\frac{B}{D}|N_k(x)|
    \leq\frac{c_0}{2}|N_k(x)|.
\]
Hence
\begin{equation}\label{eq:transition-pruned-size}
    |\M'_k|
    \geq\frac{c_0}{2}|N_k(x)|
    \qquad(k\in\mathscr K).
\end{equation}

Let
\[
    V':=\bigcup_{k\in\mathscr K}\M'_k,
    \qquad
    R:=|V'|
\]
and let $G':=G[V']$ be the subgraph of $G$ induced by $V'$.
The graph $G'$ has maximum degree at most $D$, and
\begin{equation}\label{eq:transition-total-vertices}
    R\ll_\kappa\frac{x}{y_*}.
\end{equation}

Choose a uniformly random ordering $\sigma$ of the vertices of
$G'$. For $v\in V'$, let $A_v$ denote the event that $v$
precedes every other vertex of its closed neighbourhood
$N_{G'}[v]$. Then
\[
    \mathcal I(\sigma):=\{v\in V':A_v\text{ occurs}\}
\]
is an independent set.

For $k\in\mathscr K$, put
\[
    X_k
    :=
    |\mathcal I(\sigma)\cap\M'_k|
    =
    \sum_{v\in\M'_k}\mathbf1_{A_v}.
\]
Since
\[
    \mathbb P(A_v)=\frac1{d_{G'}(v)+1}\geq\frac1{D+1},
\]
we have
\begin{equation}\label{eq:transition-expectation}
    \mu_k:=\mathbb EX_k
    \geq
    \frac{c_0}{2(D+1)}|N_k(x)|
    =:c_1|N_k(x)|.
\end{equation}

If $\sigma'$ is obtained from $\sigma$ by transposing two
vertices $u,w$, then $\mathbf1_{A_v}$ can change only if
\[
    v\in N_{G'}[u]\cup N_{G'}[w].
\]
There are at most $2(D+1)$ such vertices. Therefore
\[
    |X_k(\sigma)-X_k(\sigma')|\leq2(D+1).
\]
Theorem~\ref{ah} now gives
\begin{equation}\label{eq:transition-azuma}
    \mathbb P\left(X_k<\frac{\mu_k}{2}\right)
    \leq
    2\exp\left(
       -c_2\frac{|N_k(x)|^2}{R}
    \right)
\end{equation}
for some fixed $c_2>0$.

By \eqref{eq:transition-minimum-layer} and
\eqref{eq:transition-total-vertices},
\[
    \frac{|N_k(x)|^2}{R}
    \gg_\kappa
    \frac{x}{y_*}\,3^{-2\kappa\sqrt L+O(1)}.
\]
Since $\lambda_*<3$ for sufficiently large $x$,
\[
    y_*\leq(\log x)^{I(3)},
\]
whereas
\[
    3^{2\kappa\sqrt L}=(\log x)^{o_\kappa(1)}.
\]
Thus the exponent in \eqref{eq:transition-azuma} tends to
infinity faster than $\log L$. Since
$|\mathscr K|\ll_\kappa\sqrt L$, the union bound shows that there
exists an ordering for which
\begin{equation}\label{eq:transition-balanced-independent-set}
    |\mathcal I\cap\M'_k|
    \geq\frac{c_1}{2}|N_k(x)|
    \qquad(k\in\mathscr K).
\end{equation}

Define
\[
    F(t):=t\log 2-I(t).
\]
Then
\[
    F(2)=1,\qquad F'(2)=0,\qquad F''(2)=-\frac12.
\]
Since
\[
    \lambda_*-2=\frac{2\kappa}{\sqrt L},
\]
Taylor's theorem gives
\begin{equation}\label{eq:transition-F-expansion}
    LF(\lambda_*)
    =
    L-\kappa^2+O_\kappa(L^{-1/2}).
\end{equation}
Also,
\begin{equation}\label{eq:transition-T-expansion}
    T_*\log 2
    =
    \frac12\log L
    +O_\kappa\left(\frac{\log L}{\sqrt L}\right).
\end{equation}
It follows from \eqref{eq:transition-k0} and
\eqref{eq:transition-first-layer} that
\begin{align}
    2^{k_0}|N_{k_0}(x)|
    &\asymp_\kappa
    \frac{2^{k_0}x}{y_*} \notag\\
    &=
    x\exp\left(
       LF(\lambda_*)-T_*\log 2+O(1)
    \right) \notag\\
    &\asymp_\kappa\frac{x\log x}{\sqrt L}.
    \label{eq:transition-first-weighted-layer}
\end{align}

Furthermore,
\[
    \frac{2^{k+1}U_{k+1}}{2^kU_k}
    =\frac{2L}{k}
    =1+O_\kappa(L^{-1/2})
    \qquad(k\in\mathscr K).
\]
Since $\mathscr K$ has length $O_\kappa(\sqrt L)$,
\[
    2^kU_k\asymp_\kappa2^{k_0}U_{k_0}.
\]
Thus
\begin{equation}\label{eq:transition-each-weighted-layer}
    2^k|N_k(x)|
    \asymp_\kappa
    \frac{x\log x}{\sqrt L}
    \qquad(k\in\mathscr K).
\end{equation}

Every element of $\mathcal I\cap\M'_k$ is squarefree with exactly
$k$ prime factors. Hence, using
\eqref{eq:transition-balanced-independent-set},
\eqref{eq:transition-each-weighted-layer}, and
\eqref{eq:transition-number-layers},
\begin{align*}
    S(x,y_*)
    &\geq
    \sum_{n\in\mathcal I}\tau(n)\\
    &=
    \sum_{k\in\mathscr K}
       2^k|\mathcal I\cap\M'_k|\\
    &\gg_\kappa
    \sum_{k\in\mathscr K}2^k|N_k(x)|\\
    &\gg_\kappa x\log x.
\end{align*}
This proves \eqref{eq:transition-endpoint}.

If
\[
    1<\lambda\leq2+\frac{2\kappa}{\sqrt L},
\]
then $y_\lambda\leq y_*$ and therefore
\[
    S(x,y_\lambda)\geq S(x,y_*)\gg_\kappa x\log x.
\]
The reverse bound follows from
\[
    S(x,y_\lambda)
    \leq\sum_{n\leq x}\tau(n)
    \ll x\log x.
\]
This proves part~\textup{(ii)}.
\end{proof}

\subsection{The subcritical regime
\texorpdfstring{\(1<\lambda<2\)}{1 < lambda < 2}}
\label{subsec:subcritical}

\begin{proof}[Proof of Theorem~\ref{thm:main}\textup{(iii)}]
Fix $\lambda\in(1,2)$, and put
\[
    y:=y_\lambda=(\log x)^{I(\lambda)}.
\]
Define the fixed positive gap
\[
    \Delta_\lambda
    :=
    I(2)-I(\lambda)
    =
    (2\log 2-1)-I(\lambda)
    >0.
\]

Fix $\delta\in(0,1/2)$, and let $\varepsilon>0$. By
Lemma~\ref{lem:weighted-central}, choose fixed constants
\[
    \kappa=\kappa(\varepsilon)
    \qquad\text{and}\qquad
    M=M(\varepsilon)
\]
such that
\begin{equation}\label{eq:subcritical-central-exception}
    \sum_{\substack{n\leq x\\
                    n\notin\Nset_{\kappa,M}(x)}}\tau(n)
    \leq\varepsilon x\log x
\end{equation}
for all sufficiently large $x$. Put
\[
    V:=\Nset_{\kappa,M}(x).
\]

We first record the cardinality estimate for the corresponding
squarefree central window. Write
\[
    \Nset^{\sqf}_\kappa(x)
    :=
    \bigcup_{\substack{k\in\mathbb N\\
              |k-2L|\leq\kappa\sqrt L}}
       N_k(x),
\]
and put
\[
    k_-:=\left\lceil2L-\kappa\sqrt L\right\rceil,
    \qquad
    U_k:=\frac{x}{\log x}\frac{L^{k-1}}{(k-1)!}.
\]
For all sufficiently large $x$, every integer $k$ with
$|k-2L|\leq\kappa\sqrt L$ satisfies $L<k<3L$.
By Lemma~\ref{lem:sathe-standard} with $B=3$, we have
\[
    |N_k(x)|\ll_\kappa U_k
\]
uniformly in the central window. For all sufficiently large $x$,
$k_->3L/2$, and hence
\[
    \frac{U_{k+1}}{U_k}=\frac Lk\leq\frac23
    \qquad(k\geq k_-).
\]
It follows that
\[
    |\Nset^{\sqf}_\kappa(x)|
    \ll_\kappa U_{k_-}.
\]

Applying Robbins' inequalities
\eqref{eq:robbins-stirling} to $(k_--1)!$ gives
\[
    U_{k_-}
    \ll_\kappa
    \frac{x}{\sqrt L}
    \exp\left(
       -I\left(\frac{k_-}{L}\right)L
    \right).
\]
Since
\[
    \frac{k_-}{L}
    =
    2-\frac{\kappa}{\sqrt L}+O(L^{-1}),
\]
Taylor's theorem gives
\[
    I\left(\frac{k_-}{L}\right)
    =
    I(2)-\frac{\kappa\log 2}{\sqrt L}
    +O_\kappa(L^{-1}).
\]
Therefore
\begin{equation}\label{eq:Nsf-central-upper}
    |\Nset^{\sqf}_\kappa(x)|
    \ll_\kappa
    \frac{x\,2^{\kappa\sqrt L}}
    {(\log x)^{2\log 2-1}\sqrt L}.
\end{equation}

Let $C_3$ be the constant from Proposition~\ref{prop:good-pair-count} when $\gamma=3$, and put
\[
    \widehat C_3:=\max\{C_3,2\}.
\]
Choose $c=c_\lambda>0$ sufficiently small that
\begin{equation}\label{eq:subcritical-c-choice}
    0<c<1,
    \qquad
    c\log\widehat C_3<\frac{\Delta_\lambda}{2},
\end{equation}
and set
\[
    \ell:=\lfloor cL\rfloor.
\]
For sufficiently large $x$,
\[
    \log k\leq\ell\leq L,
    \qquad
    \ell\geq\ell_0(\delta)
\]
for every $k$ in the central window.

Define
\[
    \mathcal G_{\kappa,\delta,\ell}(x)
    :=
    \bigcup_{\substack{k\in\mathbb N\\
              |k-2L|\leq\kappa\sqrt L}}
       \M_{k,\delta,\ell}(x),
    \qquad
    \mathcal B
    :=
    \Nset^{\sqf}_\kappa(x)\setminus\mathcal G_{\kappa,\delta,\ell}(x).
\]
Choose
\[
    B_0:=I(\lambda)+1.
\]
In Proposition~\ref{prop:typical-useful-divisor}\textup{(b)},
take $\gamma=3$ and choose $A>0$ sufficiently large that
\[
    Ac>B_0+1.
\]
Since $k\asymp L$ throughout the central window, that proposition
gives
\begin{align}
    |\mathcal B|
    &\ll_{\kappa,\delta,A}
    \sqrt L\,e^{-A\ell}
    |\Nset^{\sqf}_\kappa(x)| \notag\\
    &\leq
    e^{-B_0L}|\Nset^{\sqf}_\kappa(x)|
    \label{eq:subcritical-bad-squarefree}
\end{align}
for all sufficiently large $x$.

Let $G$ be the graph on $V$ in which two distinct vertices are
adjacent when their distance is at most $y$. For
$n,n'\in V$, write
\[
    n=a\,r(n),
    \qquad
    n'=b\,r(n'),
\]
where $a,b\leq M$ are powerful and the squarefree parts
$r(n),r(n')$ lie in $\Nset^{\sqf}_\kappa(x)$.

We first count pairs for which both squarefree parts belong to
$\mathcal G_{\kappa,\delta,\ell}(x)$. For fixed $a,b$, and fixed central indices $i,j$,
Proposition~\ref{prop:good-pair-count}, applied with
\[
    K=M,\qquad \gamma=3,\qquad H=y,
\]
gives
\[
\begin{aligned}
    &\#\left\{
       (r,s)\in
       \M_{i,\delta,\ell}(x)
       \times\M_{j,\delta,\ell}(x):
       0<|ar-bs|\leq y,\ bs\leq x
    \right\}\\
    &\qquad\ll_{M,\delta}
    \frac{\widehat C_3^\ell\ell^2y}{abx}
    |N_i(x)|\,|N_j(x)|.
\end{aligned}
\]
Summing over $a,b$ and over all central $i,j$ gives
\begin{equation}\label{eq:subcritical-good-pairs}
    \#\{\text{good ordered pairs}\}
    \ll_{\kappa,\delta,M}
    \frac{\widehat C_3^\ell\ell^2y}{x}
    |\Nset^{\sqf}_\kappa(x)|^2.
\end{equation}

If at least one of the two squarefree parts lies in
$\mathcal B$, then, for fixed $a,b$ and one fixed squarefree
part, there are $O_M(y)$ possibilities for the other one.
Consequently,
\begin{equation}\label{eq:subcritical-all-pairs}
\begin{aligned}
    \#\{\text{ordered close pairs}\}
    \ll_{\kappa,\delta,M}\;&
    \frac{\widehat C_3^\ell\ell^2y}{x}
    |\Nset^{\sqf}_\kappa(x)|^2\\
    &+y|\mathcal B|.
\end{aligned}
\end{equation}
Every graph edge contributes both orientations to the count.
Thus the number of edges is bounded by the left-hand side of
\eqref{eq:subcritical-all-pairs}.

Let $P_M$ denote the number of powerful integers at most $M$.
The squarefree and bounded-powerful sets satisfy
\begin{equation}\label{eq:full-squarefree-comparable}
    |\Nset^{\sqf}_\kappa(x)|
    \leq|V|
    \leq P_M|\Nset^{\sqf}_\kappa(x)|.
\end{equation}
The lower bound comes from powerful part $1$; for the upper bound,
there are only $P_M$ possible powerful parts.

Using \eqref{eq:Nsf-central-upper},
\eqref{eq:subcritical-c-choice}, and
$\ell\leq cL$, we obtain
\begin{align*}
    \frac{\widehat C_3^\ell\ell^2y}{x}
    |\Nset^{\sqf}_\kappa(x)|
    &\ll_{\kappa,\lambda}
    L^{3/2}2^{\kappa\sqrt L}
    \exp\left(
       (c\log\widehat C_3-\Delta_\lambda)L
    \right)\\
    &\leq
    \exp\left(-\frac{\Delta_\lambda}{3}L\right)
\end{align*}
for all sufficiently large $x$. Furthermore, by
\eqref{eq:subcritical-bad-squarefree},
\[
    \frac{y|\mathcal B|}{|V|}
    \leq
    ye^{-B_0L}
    =
    \exp\bigl((I(\lambda)-B_0)L\bigr)
    =
    e^{-L}.
\]
It follows from \eqref{eq:subcritical-all-pairs} and
\eqref{eq:full-squarefree-comparable} that the number $E$ of
edges satisfies
\begin{equation}\label{eq:edge-small-iii}
    E\leq\rho_x|V|,
    \qquad
    \rho_x\leq e^{-c_\lambda' L}
\end{equation}
for some $c_\lambda'>0$ and all sufficiently large $x$.

By Theorem~\ref{cw} and Cauchy--Schwarz, there is an independent
set $\mathcal J\subseteq V$ such that
\begin{align*}
    |\mathcal J|
    &\geq
    \sum_{v\in V}\frac1{d_G(v)+1}\\
    &\geq
    \frac{|V|^2}{|V|+2E}\\
    &\geq
    (1-O(\rho_x))|V|.
\end{align*}
Thus
\[
    |V\setminus\mathcal J|\ll\rho_x|V|.
\]

If $n\in V$, its powerful part is at most $M$ and
\[
    \omega(r(n))\leq2L+\kappa\sqrt L.
\]
Consequently,
\[
    \tau(n)\ll_M2^{2L+\kappa\sqrt L}.
\]
Using \eqref{eq:Nsf-central-upper} and
\eqref{eq:full-squarefree-comparable}, we obtain
\begin{align*}
    \sum_{n\in V\setminus\mathcal J}\tau(n)
    &\ll_{\kappa,M}
    \rho_x2^{2L+\kappa\sqrt L}|V|\\
    &\ll_{\kappa,M}
    \rho_x
    \frac{x\log x}{\sqrt L}
    2^{2\kappa\sqrt L}\\
    &=o(x\log x).
\end{align*}

Finally, by \eqref{eq:subcritical-central-exception} and
\eqref{eq:dirichlet},
\begin{align*}
    \sum_{n\in\mathcal J}\tau(n)
    &\geq
    \sum_{n\leq x}\tau(n)
    -
    \sum_{\substack{n\leq x\\n\notin V}}\tau(n)
    -
    \sum_{n\in V\setminus\mathcal J}\tau(n)\\
    &\geq
    (1-\varepsilon-o(1))x\log x.
\end{align*}
Since $\mathcal J$ is independent, its distinct elements have
mutual distance greater than $y$. Hence
\[
    S(x,y)\geq(1-\varepsilon-o(1))x\log x.
\]
Because $\varepsilon>0$ is arbitrary, while
\[
    S(x,y)\leq\sum_{n\leq x}\tau(n)\sim x\log x,
\]
we conclude that
\[
    S(x,y)\sim x\log x.
\]
This proves part~\textup{(iii)}.
\end{proof}
\section*{Acknowledgements}
The author thanks his supervisor, Andrew Granville, for his invaluable guidance and suggestions throughout this work. This work was supported by a Master's Research Scholarship from the Fonds de recherche du Qu\'ebec (award no.~2004096;
\url{https://doi.org/10.69777/2004096}).

\section*{Declaration of generative AI use}
During the preparation of this manuscript, the author used ChatGPT
(OpenAI) for assistance with exposition, LaTeX formatting, literature
searches, and the review and refinement of mathematical arguments.
Claude (Anthropic) was also used to obtain comments on the draft.
The author takes full responsibility for the final manuscript,
including its mathematical content and references.


\begin{thebibliography}{99}

\bibitem{AlonSpencer}
N.~Alon and J.~H.~Spencer,
\emph{The Probabilistic Method},
4th ed., Wiley, Hoboken, NJ, 2016.

\bibitem{FordDiv}
K.~Ford,
\emph{The distribution of integers with a divisor in a given interval},
Ann. of Math. (2) \textbf{168} (2008), no.~2, 367--433.


\bibitem{Ingham}
A.~E.~Ingham,
\emph{Some asymptotic formulae in the theory of numbers},
J. London Math. Soc. \textbf{2} (1927), no.~3, 202--208.

\bibitem{Estermann}
T.~Estermann,
\emph{\"Uber die Darstellungen einer Zahl als Differenz von zwei Produkten},
J. Reine Angew. Math. \textbf{164} (1931), 173--182.

\bibitem{DI}
J.-M.~Deshouillers and H.~Iwaniec,
\emph{An additive divisor problem},
J. London Math. Soc. (2) \textbf{26} (1982), no.~1, 1--14.

\bibitem{GranvilleANT}
A.~Granville,
\emph{Analytic Number Theory Revealed: A First Guide to Prime Numbers},
Miscellaneous Books, vol.~154,
American Mathematical Society, Providence, RI, 2025.

\bibitem{FanRough}
K.~(S.)~Fan,
\emph{An inequality for the distribution of numbers free of small prime factors},
Integers \textbf{22} (2022), Paper No.~A26, 12~pp.

\bibitem{RosserSchoenfeld1962}
J.~Barkley Rosser and Lowell Schoenfeld,
\emph{Approximate formulas for some functions of prime numbers},
Illinois J. Math. \textbf{6} (1962), no.~1, 64--94.
\href{https://doi.org/10.1215/ijm/1255631807}
{\texttt{doi:10.1215/ijm/1255631807}}.

\bibitem{DudekGrytczukRucinski}
A.~Dudek, J.~Grytczuk, and A.~Ruci{\'n}ski,
\emph{Multiple twins in permutations},
Australas. J. Combin. \textbf{87} (2023), no.~3, 440--451.

\bibitem{KoukoulopoulosPrimes}
D.~Koukoulopoulos,
\emph{The Distribution of Prime Numbers},
Graduate Studies in Mathematics, vol.~203,
American Mathematical Society, Providence, RI, 2019.

\bibitem{SelbergSathe}
A.~Selberg,
\emph{Note on a paper by L.~G. Sathe},
J. Indian Math. Soc. (N.S.) \textbf{18} (1954), 83--87.
\href{https://doi.org/10.18311/jims/1954/17018}
{\texttt{doi:10.18311/jims/1954/17018}}.

\bibitem{Ten}
G.~Tenenbaum,
\emph{Introduction to Analytic and Probabilistic Number Theory},
3rd ed., Graduate Studies in Mathematics, vol.~163,
American Mathematical Society, Providence, RI, 2015.

\bibitem{Robbins}
H.~Robbins,
\emph{A remark on Stirling's formula},
Amer. Math. Monthly \textbf{62} (1955), no.~1, 26--29.
\end{thebibliography}
\end{document}